\documentclass[a4paper,10pt]{article}
\usepackage[utf8]{inputenc}
\usepackage{enumerate}
\usepackage{amssymb, bm}
\usepackage{amsthm}
\usepackage{amsmath}
\usepackage{tikz-cd}
\usepackage{tikz}
\usepackage{bbold}
\usepackage{pstricks-add, auto-pst-pdf}%
\usepackage{xcolor}

\usepackage[normalem]{ulem}
 \usepackage[all]{xy}

\newtheorem{thm}{Theorem}[section]

\newtheorem{prop}[thm]{Proposition}
\newtheorem{lemma}[thm]{Lemma}
\newtheorem{remark}[thm]{Remark}

\newtheorem{defn}[thm]{Definition}

\newcommand{\re}[1]{(\ref{#1})}

\newcommand{\rl}[1]{Lemma~\ref{#1}}

\newcommand{\rp}[1]{Proposition~\ref{#1}}

\newcommand{\rd}[1]{Definition~\ref{#1}}
\newcommand{\rrem}[1]{Remark~\ref{#1}}
\newcommand{\diag}{\operatorname{diag}}

\newcommand{\Om}{\Omega}

\def\Z{\mathbb{Z}}
  \def\C{\mathbb{C}}
\def\R{\mathbb{R}}
\def\N{\mathbb{N}}

\def\d{\partial}

\def\cA{\mathcal{A}}
  \def\cO{\mathcal{O}}

\def \Im{\mathrm{Im}}
\def \Re{\mathrm{Re}}

\newcommand{\dindice}[2]{{\stackrel{\scriptstyle #1}{\scriptstyle #2}}}

\title{Ueda foliation problem for complex tori}
\author{Laurent Stolovitch\thanks{CNRS and Laboratoire J.-A. Dieudonn\'e
		U.M.R. 7351, Universit\'e Côte d'Azur, Parc Valrose 06108 Nice Cedex 02, France, email: {\tt stolo@unice.fr}.} $\;$ and$\;$  Xiaojun Wu\thanks{Laboratoire J.-A. Dieudonn\'e,	U.M.R. 7351, Universit\'e Côte d'Azur, Parc Valrose	06108 Nice Cedex 02, France, email: {\tt Xiaojun.WU@univ-cotedazur.fr}. This work was supported by the French government through the France 2030 investment plan managed by the National Research Agency (ANR), as part of the Initiative of Excellence of Université Côte d’Azur under reference number ANR-15-IDEX-01.} }

\begin{document}
\maketitle
\begin{abstract}
	We consider an embedded general complex torus $C_n$ into a complex manifold $M_{n+d}$ with a unitary flat normal bundle $N_C$. We show the existence of (non-singular) holomorphic foliation in a neighborhood of $C$ in $M$ having $C$ as leaf under some conditions.
\end{abstract}
\section{Introduction}
Let $\iota: C \to M $ be a holomorphic embedding of a compact complex
manifold $C$ of dimension $n$ in a complex manifold $M$ of dimension $n+d$. We shall still denote $\iota (C)$ by $C$. Let $N_C$ be its normal bundle. We assume that $TM_{|C}$ splits, that is $TM_{|C}=N_C\oplus TC$ and that $N_C$ is unitary flat, that is admit locally constant unitary transition matrices. We aim at giving sufficient conditions ensuring the existence of a holomorphic foliation having $C$ as a leaf in some neighborhood of $C$ in $M$. T. Ueda \cite{Ued82} studied the case of an embedded complex compact curve into a surface and showed, in the so called {\it infinite type case}, the existence of such a foliation under a ``Diophantine-like" condition of the form: there exist $M>0$,$\tau\geq 0$ such that for all $l>2$,  $\text{dist}(\mathbb{1}, N_C^{-l+1})>Ml^{-\tau}$. Here the distance is the one defined on the Picard group of $C$. Recently,
the problem of existence of holomorphic foliation in a neighborhood of an embedded compact manifold $C$ of which it is a leaf has attracted lot of attention (e.g. \cite{MR1967036,koike-ueda0,koike-ueda, CLPT,GS21}). 
The aim of this article is to provide with a new range of such examples, namely embedded general tori in complex manifold of any dimension.
This problem is also related another one : the existence of a neighborhood biholomorphic to a neighborhood of the zero section in its normal bundle. The latter is related to Grauert's ``Formale Prinzip" in the case of a flat normal bundle. This situation is quite different than when there is ``curvature",  as initiated by Grauert \cite{MR0137127} (negative case) or Griffiths \cite{MR206980} (positive case). This was first devised by V.I Arnold\cite{arnold-embed} for elliptic curve embedded in surfaces and generalized to an abstract situation by the first author and X. Gong \cite{GS21}. Very recently  the problem of equivalence of neighborhoods was solved by the first author and X. Gong for embedded general complex tori \cite{GS}.

The main result of this article is the following.
\begin{thm}
	Let $C$ be an $n$-dimensional complex torus, holomorphically embedded into a complex manifold $M_{n+d}$.   Assume that $T_M|_C$ splits. Assume the normal bundle $N_C$ has (locally constant) unitary transition functions. Assume that  $N_C$ is vertically  strongly Diophantine (see Definition (\ref{dioph})). Then there exists a non-singular
	holomorphic foliation of the germ of neighborhood of $C$ in $M$ having $C$
	as a compact leaf.
\end{thm}

In this note, following the approach of \cite{GS}, we will show that under some arthimetic assumption on some Stein covering of a complex torus, we can ``vertically linearize" a (holomorphic) neighborhood of the torus (and thus show the existence of  a nonsingular holomorphic foliation which have the torus as a compact leaf).

The principal technical novelty is the following. Rather than employing coverings by finite open sets and cocycle-type arguments, as inspired by Ueda, to attain an $L^\infty$ estimate on a larger domain that is independent of the recurrence procedure (see \cite[p.37, (3.15)]{GS21}), we apply a Hartogs-type lemma to the double translations in $n$ directions of the lattice of the fundamental domain. The selection of  these translations facilitates the embedding of the translation of a larger domain (than the fundamental domain) within the convex hull of the  union of all these translated domains.

The organisation of the article is as follows.
 In Section 2, we revisit the problem's context as outlined in \cite{GS}, and we establish essential lemmas for subsequent use. Notably, the Hartogs-type lemmas \ref{convh} and \ref{cauchy-transl} are elementary, yet they appear to be novel in the context of the linearization problem, capitalizing on the distinct features of the complex torus. Moving on to Section 3, we provide a detailed proof of the main result.
\\
{\bf Acknowledgment~:} The first author thanks X. Gong for discussions on the subject.
\section{Setting}
Let $U$ be a neighborhood of $C$ in $M$ such that $U$ admits a smooth, possibly non-holomorphic strong retract to $C$; namely there is a smooth mapping $R\colon U\times [0,1]\to U$ such that $R(\cdot,0)= \mathrm{Id}$ on $U$, $R(\cdot,t)= \mathrm{Id}$ on $C$, and $R(\cdot,1)(U)=C$. Thus, $\pi_1(U,x_0)=\pi_1(C,x_0)$ for $x_0\in C$.
Since
we are considering only unprecised neighborhoods of $C$ in $M$, we will identify $U$ and $M$.
Recall the following lemma in \cite[Lemma 4.1]{GS} which relate the covering of the submanifold and the covering of its neighborhood.
\begin{lemma}
Let $C$ be a compact complex manifold. Let $\pi\colon \tilde{C}\to C$ be a holomorphic covering and $\pi(x_0^*)=x_0$. Suppose that $(M,C)$ is a holomorphic neighborhood of $C$. There is a neighborhood $U$ in $M$ of $C$ 
 and  a holomorphic neighborhood  $\tilde U$ of $\tilde{C}$ such that $p\colon\tilde U\to U$  is an extended covering of the covering $\pi\colon \tilde{C}\to C$ and $C$
(resp. $\tilde{C}$)  is a smooth strong retract of $U$ (resp. $\tilde U$). Consequently, $$
  \pi_1(\tilde U,x_0^*)=\pi_1(\tilde{C},x_0^*), \quad \pi_1(U,x_0)=\pi_1(C,x_0).
$$
\end{lemma}
Applying the above lemma to $(N_C,C)$ and a covering $\pi|_{\tilde{C}}\colon \tilde{C}\to C$, we have a covering $\hat\pi\colon \widetilde{N_C}\to N_C$ such that
$$
\tilde{C}\subset\widetilde{N_C}, \quad \pi_1(\widetilde{N_C},x_0^*)=\pi_1(\tilde{C},x_0^*), \quad \pi_1(N_C,x_0)=\pi_1(C,x_0).
$$

In the following, we will consider the case that $C$ is a complex torus.
We will always add this assumption from now on.

Consider the fundamental domain
$$
\omega_0=\left\{\sum_{j=1}^{2n} t_je_j\in\C^n\colon t\in[0,1[^{2n}\right\},
$$
where $e_i(1 \leq i \leq 2n)$ is chosen such that $C \simeq \C^n / \Lambda$ where $\Lambda:= \sum_{1 \leq i \leq 2n} \Z e_i$.
Without loss of generality, we can assume that $e_i(1 \leq i \leq n)$ is the standard base of $\C^n$ (as a complex vector space).

Consider the cylinder
$$\tilde{C}:= \C^n/ \sum_{1 \leq i \leq n} \Z e_i$$
such that the natural quotient map gives $\pi\colon \tilde{C}\to C$ a holomorphic covering.
Note that by \cite[Proposition 3.6]{GS}, $\widetilde{N_C}$ is (holomorphically) trivial if $N_C$ is flat.

For $\epsilon>0$,  define the Reinhardt domain ${\Om 
}_\epsilon$ by
\begin{align*}
\omega_\epsilon&:=\left\{\sum_{j=1}^{2n} t_je_j\colon t=(t_1,\ldots,t_{2n})\in[0,1[^n\times ]-\epsilon,1+\epsilon[^{n}\right\},\\
{\Om}_\epsilon&:=\{(e^{2\pi i\zeta_1}, \dots e^{2\pi i\zeta_n})\colon\zeta\in \omega_\epsilon\}.\\
 \Omega_{\epsilon}^+& =\{(|z_1|,\ldots,|z_n|), z\in \Omega_{\epsilon}\}\\
&= \left\{(e^{-2\pi R_1},\ldots,e^{-2\pi R_n}),\,R=\sum_{i=1}^nt_{i+n}\Im\, e_{i+n},\, t''\in]-\epsilon,1+\epsilon[^n\right\}.
\end{align*}
With $\Delta_r=\{z\in\C\colon|z|<r\}$, we also define 
\begin{equation}\label{domains}
\omega_{\epsilon,r}:= \omega_\epsilon\times  \Delta^d_r, 
\quad \Omega_{\epsilon,r}:=\Omega_\epsilon\times  \Delta^d_r. 
\end{equation}
 A function on $\omega_{\epsilon,r}$ that has period $1$ in all $z_j$ is identified with a function on $\Omega_{\epsilon,r}$.
 In the following,  we will denote by $(h,v)$ the coordinates of $\Omega_{\epsilon,r}$.
Also we will call the $h$-components (resp. $v$-component) of an element of $\Omega_{ \epsilon,r}$, its horizontal (resp. vertical) component.

 In the case of torus, we have the following result \cite[Proposition 4.3]{GS} on the classification of pair $(C,M)$.
\begin{prop}
Let $C$ be the complex torus and $\pi \colon\tilde{C}=\C^n/{\Z^n}\to C$ be the covering defined above.
Let $(M,C)$ be a neighborhood of $C$. Assume that $N_C$ is flat.
Then
$(M,C)$ is holomorphically equivalent to the quotient space of an open neighborhood of $\tilde{C}$ in $\widetilde{N_C}$  by the deck transformations of $\widetilde{N_C}$.
Moreover,
one can take $\omega_{\epsilon_0,r_0}$ (for suitable choice of $\epsilon_0,r_0$) such that $(M,C)$ is biholomorphic to the quotient of $\omega_{\epsilon_0,r_0}$ by  $\tau^0_1,\dots,\tau^0_n$.
Here for any $1 \leq i \leq n$, $\tau^0_i$ is the translation by $e_{i+n}$ when restricted to 
$\omega_\epsilon \times \{0\}$.
Let $\tau_j$ be the mapping defined on $\Omega_{\epsilon_0,r_0}$ corresponding to $\tau_j^0$.
Then $\tau_1,\dots, \tau_n$ commute pairwise   wherever they are defined, i.e.
$$
\tau_i\tau_j(h,v)=\tau_j\tau_i(h,v)\quad  \forall  i\neq j
$$
for $(h,v)\in \Omega_{\epsilon_0, r_0}\cap \tau_i^{-1}\Omega_{\epsilon_0,r_0}\cap \tau_j^{-1}\Omega_{\epsilon_0,r_0}$.
Notice that
$(M,C)$ is also biholomorphic to the quotient of $\Omega_{\epsilon_0,r_0}$ by  $\tau_1,\dots,\tau_n$.

On the other hand, let $(\tilde M,C)$ be another  such neighborhood having the corresponding generators $\tilde\tau_1,\dots,\tilde\tau_n$  defined on $\Omega_{\tilde\epsilon_0,\tilde r_0}$. Then $(M,C)$ and $(\tilde M,C)$ are holomorphically equivalent if and only if there is a  biholomorphic mapping $F$  from $\Omega_{\epsilon,r}$ into $\Omega_{\tilde \epsilon,\tilde r}$ for some positive $\epsilon,r,\tilde \epsilon,\tilde r>0$  such that
$$
F \tilde\tau_j(h,v)=\tau_jF(h,v),\ j=1,\ldots, n,
$$
wherever   both sides are defined, i.e. $(h,v)\in \Omega_{\tilde\epsilon,\tilde r}\cap \tilde\tau_j^{-1}\Omega_{\epsilon,r}\cap \Omega_{\epsilon,r}\cap F^{-1}\Omega_{\epsilon,r}.$
\end{prop}
 By identification of $\C^n/\Z^n=(\C/ \Z)^n=(\C^*)^n$,
 we can identify $\widetilde{N_C}$ as $(\C^*)^n \times \C^d$ (since $\widetilde{N_C}$ is holomorphically trivial).
 We consider $\Omega_{\epsilon, r}$ as an open subset of $(\C^*)^n \times \C^d$.
 We assume from now on that $N_C$ is Hermitian flat.
 The deck transformations of $\widetilde{N_C}$ in this case (cf. \cite[(4.7), (4.8)]{GS}) can be chosen to be given by for any $1 \leq j \leq n$,
\begin{equation} \hat \tau_j(
 	h,v)
 	=
 	(T_jh,  M_jv) \label{tauhat},\quad\hat\tau_0(h,v):=(h,v)
\end{equation}
for some diagonal matrix
 $$
 T_j:= \diag(\lambda_{j,1},\dots, \lambda_{j,n}), M_j:=\diag(\mu_{j,1},\dots, \mu_{j,d})$$
and $(h,v) \in (\C^*)^n \times \C^d.$
We assume from now on that $T_M|_C$ splits (i.e. $T_M|_C=T_C \oplus N_C$).
By construction, for $(h,v)\in \Omega_{\epsilon, r}$,
$$
      \tau_j(h,0)=( T_jh,0).
      $$
    Since $T_M|_C$ splits, the differential of $\tau_j$ along $\tilde{C}$
give the deck transformations of  $\widetilde{N_C}$.
In other words, $\tau_j$  in the horizontal direction is a higher-order perturbation of 
      $(T_j,M_j)$ (of order $\geq 2$ in $v$).

Recall the notations
$
\lambda_{l}=(\lambda_{l,1},\ldots,\lambda_{l,n})$ and $ \mu_l=(\mu_{l,1},\ldots,\mu_{l,d}).
$
With $P=(p_1,\ldots, p_n)\in\Z^n$ and $Q=(q_1,\ldots, q_d)\in \mathbb{N}^d$, we define
$$
\lambda_{l}^P\mu_{l}^Q:=\prod_{i=1}^n\lambda_{l,i}^{p_i}\prod_{j=1}^d\mu_{l,j}^{q_j}.
$$
We will need the following sufficient condition to vertically linearize the Deck transformations.
\begin{defn}\label{dioph}
The pullback 
normal bundle $\widetilde N_C$  is said to be  {vertically Diophantine} (resp. strongly Diophantine) if for all $(Q,P)\in \mathbb{N}^d\times \Z^n$,   $|Q|>1$ and all   $j=1,\ldots ,d$,
\begin{equation}
\max_{l\in \{1,\ldots, n\}}  \left |\lambda_l^P 
\mu_{l}^Q-\mu_{{l},j}\right |   >  \frac{D}{(|P|+|Q|)^{\tau}},\label{dv}
\end{equation}
(resp.  
\begin{equation}
\forall 1\leq l \leq n,\quad	\left |\lambda_l^P 
	\mu_{l}^Q-\mu_{{l},j}\right |   >  \frac{D}{(|P|+|Q|)^{\tau}}\label{dv-str}
\end{equation})
for some $D>0$, $\tau>0$ (independent of $P,Q$).
\end{defn}
\begin{remark}
The vertically (resp. strongly) Diophantine condition is independent of choice of generators as shown in \cite[Proposition 4.7]{GS}.
In particular, it is equivalent to the condition that
for all $(Q,P)\in \mathbb{N}^d\times \Z^n$,   $|Q|>1$ and all   $j=1,\ldots ,d$,
\begin{equation}
\max_{l\in \{1,\ldots, n\}}  \left |\lambda_l^{-P} 
\mu_{l}^{-Q}-\mu_{{l},j}^{-1}\right |   >  \frac{D}{(|P|+|Q|)^{\tau}}\label{dv-}
\end{equation}
(resp.
\begin{equation}
	\forall l=1,\ldots n,\quad  \left |\lambda_l^{-P} 
	\mu_{l}^{-Q}-\mu_{{l},j}^{-1}\right |   >  \frac{D}{(|P|+|Q|)^{\tau}}\label{dv-}
\end{equation})
for some $D>0$, $\tau>0$ (independent of $P,Q$).

If $\widetilde N_C$ is vertically strongly Diophantine, then it is also strongly non-resonnant~:
$$
\forall (Q,P)\in \mathbb{N}^d\times \Z^n, |Q|>1, j=1,\ldots ,d, l=1,\ldots, n,\quad \lambda_l^P 
\mu_{l}^Q-\mu_{{l},j}\neq 0.
$$
\end{remark}
\begin{remark}
Note that without any arthimetic assumption, Ueda \cite{Ued82} showed an example where there exists no regular foliation near an elliptic curve $C$ in a complex surface such that the elliptic curve is a leaf of the foliation, the tangent bundle of the surface splits and the normal bundle of the elliptic curve is trivial.
In this example, there exist compact irreducible curves $C_{i}$ such that $C_i$ is cohomolgous to $m_i C$ with $\lim_{i \to \infty} m_i=\infty$.
However, if the claimed regular foliation exists, the leaves would be the only irreducible compact curves near $C$.
\end{remark}
\begin{defn}\label{def-domains}
Set
 $\Omega_{\epsilon,r}:=\Omega_\epsilon\times\Delta_r^d$
 ,
and for $\ell\in\N$ and $1\leq i\leq n$~:
\begin{align}
	\tilde\Om^{(\pm \ell)}_{i,\epsilon,r}&:=\Om_{\epsilon,r}\cup\hat\tau_i^{\pm 1}(\Om_{\epsilon,r})\cup \cdots\cup \hat\tau_i^{\pm \ell}(\Om_{\epsilon,r}),\nonumber\\
	\tilde\Om^{(\pm \ell)}_{\epsilon,r}&=\cup_{i=1}^n\tilde\Om^{(\pm \ell)}_{i,\epsilon,r}\label{domainsl}.
\end{align}
Denote by $\cA_{\epsilon,r}$
(resp. $\tilde{\cA}^{(\pm\ell)}_{\epsilon, r}$)
the set of holomorphic functions on a neighborhood of $\Omega_{\epsilon,r}$,
(resp. $\tilde \Omega^{(\pm\ell)}_{\epsilon,r}$).
If $f\in \mathcal{A}_{\epsilon, r}$,
we set
$$
||f||_{\epsilon,r}:=\sup_{(h,v)\in\Omega_{\epsilon,r}}|  f(h,v)|.
$$
More generally, if $U$ is open subset in $(\C^{*})^n \times \C^d$ and $f$ is holomorphic in a neighborhood of $U$, then $$\|f\|_U:=\sup_{(h,v)\in U}|  f(h,v)|.$$
\end{defn}
As such, each
  $f\in \mathcal{A}_{\epsilon, r}$
can be expressed as a convergent Taylor-Laurent series
$$
f(h,v)= \sum_{P\in \Z^n, Q \in \N^d}f_{Q,P}h^Pv^Q
$$
for $(h,v)\in \Omega_{\epsilon,r}$.

We have the following basic Cauchy estimate.
\begin{lemma}
\label{cauchy}
 If $f=\sum_{P\in \Z^n}f_{P}(v)h^P=\sum_{Q\in \N^d}f_{Q}(h)v^Q \in \cA_{\epsilon , r}$   
 and  $0<\delta<{\kappa}\epsilon$, $\delta \leq \delta_0$,
	then
	\begin{equation}
\label{Reps}
	|f_{Q}(h)|\leq \frac{\sup_{\Omega_{\epsilon,r}}|f|}{r^{|Q|}}.
\end{equation}
	\begin{equation}
	\label{Reps2}
	\|f\|_{\epsilon
  {-{\delta}/{\kappa}}, r}\leq \frac{C\sup_{\Omega_{\epsilon,r}}|f|}{\delta^{\nu}},
	\end{equation}
		\begin{equation}
	\label{Reps3}
	\|\partial^{P_0}_h f\|_{\epsilon
  {-{\delta}/{\kappa}}, r}\leq \frac{CC'^{|P_0|} \sup_{\Omega_{\epsilon,r}}|f|}{\delta^{\nu+|P_0|}},
	\end{equation}
	where $C$ and $\nu$ depends only on $n$ and $\delta_0$.
Here, $\kappa$ is some constant independent of $\epsilon,r$ and $C'$ depends only on $ \epsilon$.
\end{lemma}
\begin{proof}
The proof of the estimates $(\ref{Reps})$ is given in Lemma 4.15 of \cite{GS}.
To get the estimate $(\ref{Reps2})$  and $(\ref{Reps3})$,
we can modify the proof of Lemma 4.15 of \cite{GS} as follows.
We recall the notation
$$
{\mathcal P}_\epsilon^+=\left\{\sum_{i=1}^n t_i\Im \tau_i\in\R^n\colon t\in ]-\epsilon,1+\epsilon[^{n}\right\},$$$$
{\Om}_\epsilon^+=\left\{(e^{-2\pi R_1},\dots, e^{-2\pi R_n})\colon R
\in {\mathcal P}_\epsilon^+
\right\}.$$
According to \cite[Lemma 4.12]{GS} and Cauchy estimates for polydiscs, we have if $(h,v)\in \Omega_{\epsilon -\delta,r}$, then for all $s\in \Omega_{\epsilon}^+$ and any fixed $v$,
\begin{equation}\label{coeffhP}
	|f_P(v)h^P|\leq \left|\frac{1}{(2\pi i)^n}\int_{|\zeta_1|=s_1,\dots, |\zeta_n|=s_n} f(\zeta,v)\frac{h^P}{\zeta^{P}}\, \frac{d\zeta_1\wedge\cdots \wedge d\zeta_n}{\zeta_1\cdots\zeta_n}\right|,
\end{equation}
Set $s_j=e^{-2\pi R_j}$, $|h_j|=e^{-2\pi R_j'}$, $R=(R_1,\cdots, R_n)$ and $R'=(R_1',\cdots, R_n')$.  By \cite[Lemma 4.13]{GS},
		\begin{equation}
\label{fourier-type}
\inf_{(|\zeta_1|,\cdots,|\zeta_n|)=s\in\Omega_{\epsilon}^+}\sup_{h\in \Omega_{\epsilon  {-\delta}}}
\left|\frac{h^P}{\zeta^P}\right|=\inf_{R\in \mathcal P_{\epsilon}^+}\sup_{R'\in\mathcal P_{\epsilon  {-\delta}}^+} e^{-2\pi <R-R',P>}\leq e^{-\kappa\delta' |P|},
		\end{equation}
where the positive constant $\kappa$ depends only on $\Im \tau_i$
and $\delta'=\delta/ \kappa$.  Thus
		\begin{equation}
\label{L^inf}
| f_{P}(v)h^P|\leq \sup_{\Omega_{\epsilon,r}}|f| e^{-\delta |P|}.
\end{equation}
Similarly, we have
$$
|\partial^{P_0}_h f(h,v)|\leq\sum_{P \in \Z^n}  \left|\binom{P}{P_0} f_{P}(v)h^{P-P_0}\right|
$$
where $P_0=(P_{0,1}, \cdots, P_{0,n})$ and
$$\binom{P}{P_0}= \prod_{j=0}^n \prod_{i=0}^{P_{0,j}-1} (P_j-i).$$
The estimate follows by summing and using (\ref{L^inf}) which gives
$$
|\partial^{P_0}_h f(h,v)|\leq C  \sup_{\Omega_{\epsilon,r}}|f|  \prod_{i=1}^n s_i^{P_{0,i}} /\delta^{\nu+|P_0|}.
$$
\end{proof}
 \begin{remark}\label{rem-dom}
Assume that each $M_k$ is a unitary matrix. Then, all estimates of \rl{cauchy} remain valid if one replaces $\Om_{\epsilon,r}$ by $\hat\tau_k(\Om_{\epsilon,r})$ for any $1\leq k\leq n$.
\end{remark}
We will need the following type of Hartogs lemma.
\begin{lemma}\label{convh}
Fix $\epsilon >0$.
Define
$$\Omega'_\epsilon:= \Omega_\epsilon \cup \cup_i (T_i \Omega_\epsilon \cup T^2_i \Omega_\epsilon) \cup \cup_i( T_i^{-1} \Omega_\epsilon \cup T_i^{-2} \Omega_\epsilon)$$
which is a Reinhardt domain (i.e. a domain such that $ (e^{\sqrt{-1} \theta_1} z_1, \cdots , e^{\sqrt{-1} \theta_n} z_n)$ is in $\Omega'_\epsilon$ for
every $z = (z_1, \cdots , z_n) \in \Omega'_\epsilon$ and $\theta_1, \cdots \theta_n \in \R$).
Consider its
logarithmic indicatrix
$\omega'^*_\epsilon=\Omega'^*_\epsilon \cap \R^n$
with
$\Omega'^*_\epsilon=\{\xi \in \C^n;(e^{\xi_1}, \cdots, e^{\xi_n}) \in \Omega'_\epsilon\}$.
Denote $\varphi(z_1, \cdots , z_n):=(\log| z_1|, \cdots ,\log| z_n|)$.
Let $F \in \cO(\Omega'_\epsilon)$.
The $F$ can be extended over the preimage of the convex hull $\mathrm{Conv}(\omega'^*_\epsilon)$ of $\omega'^*_\epsilon$ under $\varphi$.
Moreover, the $L^\infty$ norm of extended function is equal to the $L^\infty$ norm of $F$ on $\Omega'_\epsilon$.
\end{lemma}
\begin{proof}
Consider the following mapping $\psi(z_1, \cdots , z_n)=(e^{z_1},  \cdots ,e^{ z_n})$.  We have for $z\in\C^n$, $\theta\in\R^n$
$$
\varphi\circ\psi(z)=(\Re\, z_1,\ldots, \Re\, z_n),\quad \varphi(\zeta_1e^{i\theta_1},\ldots, \zeta_ne^{i\theta_n})=\varphi(\zeta).
$$
Since $\Omega'_\epsilon$ is a Reinhardt domain, if $\zeta=\psi(z)\in \Omega'_\epsilon$, then for $\theta\in\R^n$, $$(e^{z_1}e^{i\theta_1},\ldots,e^{z_n}e^{i\theta_n})=\psi(z+i\theta)\in \Omega'_\epsilon.$$ Hence, $\Omega'^*_\epsilon$ is a tube and we have
\begin{align*}
	\psi^{-1} \varphi^{-1}(\omega'^*_\epsilon)&=\omega'^*_\epsilon+ \sqrt{-1}\R^n=\Omega'^*_\epsilon.\\
	 \psi(\Omega'^*_\epsilon)&= \Omega'_\epsilon,\quad \varphi(\Omega'_\epsilon)=\omega'^*_\epsilon.
\end{align*}
Let us set $\tilde \Omega'_{\epsilon}:=\psi(\mathrm{Conv}(\omega'^*_\epsilon)+\sqrt{-1} \R^n) $.
We can summerize the inclusion of open sets as follows.
$$\begin{tikzcd}
\Omega'^*_\epsilon=\omega'^*_\epsilon+\sqrt{-1} \R^n \arrow[d, "\mathrm{etale}"] \arrow[r, "\subset"] & \mathrm{Conv}(\omega'^*_\epsilon)+\sqrt{-1} \R^n \arrow[r, "\subset"] \arrow[d, "\mathrm{etale}"] & \C^n \arrow[d, "\psi"]    \\
 \Omega'_\epsilon \arrow[d] \arrow[r, "\subset"]                      & \tilde\Omega'_\epsilon \arrow[r, "\subset"] \arrow[d]                      & \C^n \arrow[d, "\varphi"] \\
\omega'^*_\epsilon \arrow[r,"\subset"]                                             
 & \mathrm{Conv}(\omega'^*_\epsilon) \arrow[r, "\subset"]                                              & \R^n
\end{tikzcd}$$
Let us consider the function $\psi^* F$. It is defined on $\Omega'^*_{\epsilon}$ which is the tube over $\omega'^*_{\epsilon}$. By a result of Bochner \cite{Bo38}, it extends over the tube over the convex hull $\mathrm{Conv}(\omega'^*_\epsilon)$  of $\omega'^*_\epsilon$. It is also the preimage of the convex hull $\mathrm{Conv}(\omega'^*_\epsilon)$ under $\varphi \circ \psi$.
On the other hand, for each variable, $\psi^* F$ is $2 \pi \sqrt{-1}$ periodic over  the preimage of  $\omega'^*_\epsilon$ under $\varphi \circ \psi $ which is an open set in the connected preimage of its convex hull.
By identity theorem, the extension is unique and for each variable, $\psi^* F$ is $2 \pi \sqrt{-1}$ periodic
which defines a holomorphic function on  the preimage of the convex hull $\mathrm{Conv}(\omega'^*_\epsilon)$ of $\omega'^*_\epsilon$ under $\varphi$. Hence, $F$ extends to holomorphic function on $\tilde \Omega'_{\epsilon}$.

For the last statement,
note that $\psi^* |F|^2$ is subharmonic over $\omega'^*_\epsilon+ \sqrt{-1}\R^n$.
For any point in $\mathrm{Conv}(\omega'^*_\epsilon)$, by definition, there exist $x, y \in \omega'^*_\epsilon$ such that the segment $L$ connecting $x, y$ is contained in $\mathrm{Conv}(\omega'^*_\epsilon)$ containing the given point (see Figure \ref{figure-convh} below).
The restriction of $\psi^* |F|^2$ over $L+\sqrt{-1}\R^n$ is still subharmonic.
In particular, by mean value inequality involving Poisson kernel (see e.g. \cite[(4.12), Chap. I]{agbook}),
$$\sup_{L+\sqrt{-1}\R^n} \psi^* |F|^2=\sup_{\d L+\sqrt{-1}\R^n} \psi^* |F|^2$$
which implies
$$\sup_{\omega'^*_\epsilon+\sqrt{-1}\R^n} \psi^* |F|^2=\sup_{\mathrm{Conv}(\omega'^*_\epsilon)+\sqrt{-1}\R^n} \psi^* |F|^2.$$
This finishes the proof of last statement.
\end{proof}
\begin{lemma}\label{cauchy-transl}
Under the same notations of previous lemma \ref{convh}, there exists $\eta>0$ depending on $\epsilon$ such that
$$ \cup_i T_i \Omega_{\epsilon+\eta} \cup \cup_i T_i^{-1} \Omega_{\epsilon+\eta}$$
is contained in the preimage of the convex hull $\mathrm{Conv}(\omega'^*_\epsilon)$ of $\omega'^*_\epsilon$ under $\varphi$.
\end{lemma}
\begin{proof}
The statement is equivalent to $$\varphi(\cup_i T_i \Omega_{\epsilon+\eta} \cup \cup_i T_i^{-1} \Omega_{\epsilon+\eta}) \subset \mathrm{Conv}(\omega'^*_\epsilon)$$
for some $\eta >0$,
which is
invariant under the base change of $\R^n$.
In particular, without loss of generality, we may assume that $T_j$ corresponds to translation  $\tilde T_j$ by $e_j$ the standard basis of $\R^n$ and
$$\omega_\epsilon=]-\epsilon,1+\epsilon[^n$$
where it is easy to check the statement.

The corresponding picture in dimension 2 is as follows.
The domain bounded by dashed lines is $ \mathrm{Conv}(\omega'^*_\epsilon)$.
The greyed domain
is $\tilde T_1 \omega^*_{\epsilon+\eta}$.
The union of domains bounded by solid lines is
$ \omega'^*_\epsilon= \omega^*_\epsilon \cup \cup_i (\tilde T_i \omega^*_\epsilon \cup \tilde T^2_i \omega^*_\epsilon) \cup \cup_i( \tilde T_i^{-1} \omega^*_\epsilon \cup \tilde T_i^{-2} \omega^*_\epsilon)$.
To ease the illustration, we take $\epsilon=0$ in the picture.
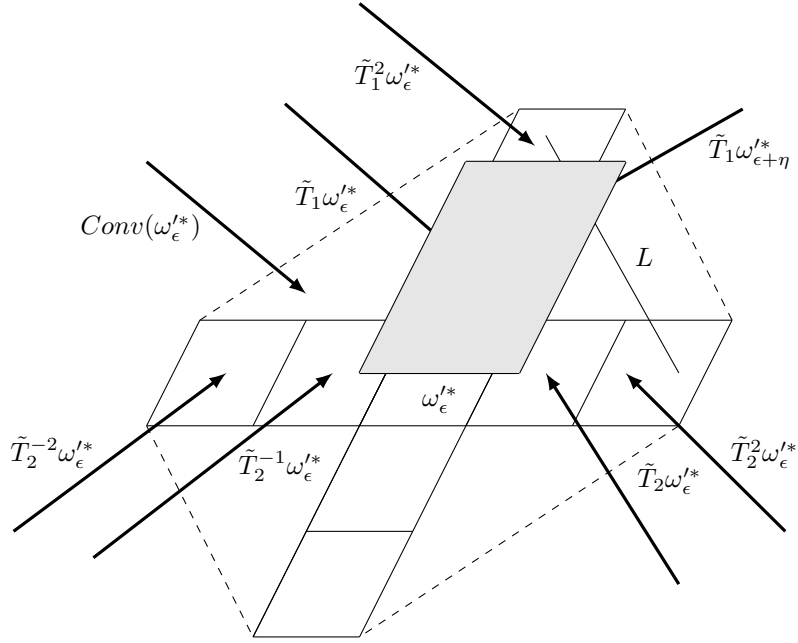
\begin{figure}
 \begin{tikzpicture}[scale=.7]
 	\tikzstyle{estun}=[->,>=latex,very thick,dotted]
 	\tikzstyle{estunerc}=[->,>=latex,very thick,dotted,olive]
 	\tikzstyle{implique}=[->,>=latex,very thick]
 	\tikzstyle{impliquel}=[<-,>=latex,very thick]
 	\tikzstyle{equiv}=[<->,>=latex,very thick, dotted]
 	\tikzstyle{equivlr}=[<->,>=latex,very thick]
\draw (0,0) -- (10,0) ;
\draw (0,0) -- (1,2) ;
\draw (1,2) -- (11,2) ;
\draw (10,0) -- (11,2) ;
\draw (2,0) -- (3,2) ;
\draw (4,0) -- (5,2) ;
\draw (6,0) -- (7,2) ;
\draw (8,0) -- (9,2) ;
\draw (2,-4) -- (7,6) ;
\draw (4,-4) -- (9,6) ;
\draw (6,4) -- (8,4) ;
\draw (7,6) -- (9,6);
\draw (2,-4) -- (6,4) ;
\draw (3,-2) -- (5,-2) ;
\draw (2,-4) -- (4,-4) ;
\draw[dashed] (1,2) -- (7,6) ;
\draw[dashed] (11,2) -- (9,6) ;
\draw (10,1) -- (7.5,5.5) ;
\draw (9,3.2) node[right]{$L$};
\draw[implique](10,-3) -- (7.5,1)node[midway,right]{$\;\;\tilde T_2 \omega'^*_{\epsilon}$};
\draw[implique](12,-2) -- (9,1)node[midway,right]{$\;\;\tilde T_2^2 \omega'^*_{\epsilon}$};
\draw[implique](-1,-2.5) -- (3.5,1)node[midway,right]{$\;\;\tilde T_2^{-1} \omega'^*_{\epsilon}$};
\draw[implique](-2.5,-2) -- (1.5,1)node[midway,left]{$\tilde T_2^{-2} \omega'^*_{\epsilon}\;\;$};
\draw[implique](2.6,6.1) -- (6.6,2.6)node[midway,left]{$\tilde T_1 \omega'^*_{\epsilon}\;\;\;$};
\draw[implique](4,8) -- (7.3,5.3)node[midway,left]{$\tilde T_1^2 \omega'^*_{\epsilon}\;\;$};
\draw[implique](11.2,6) -- (8.4,4.4)node[midway,right]{$\;\;\;\;\tilde T_1 \omega'^*_{\epsilon+\eta}$};
\draw[implique](0,5) -- (3,2.5)node[midway,left]{$Conv(\omega'^*_{\epsilon})\;\;$};
\draw[dashed] (0,0) -- (2,-4) ;
\draw[dashed] (10,0) -- (4,-4) ;
\draw[black, thick] (4,1) -- (7,1) ;
\draw[black, thick] (6,5) -- (9,5) ;
\draw[black, thick] (4,1) -- (6,5) ;
\draw[black, thick] (7,1) -- (9,5) ;
\draw (5,0.5) node[right]{$\omega'^*_\epsilon$};
\draw [black, thick](5.7,3.3);
\fill[color=gray!20,opacity=0.5]
(4,1) -- (7,1)  -- (9,5)-- (6,5)--cycle;
\end{tikzpicture}
\caption{Domains and their translates by the $\tilde T_j$'s for $\epsilon=0$ and their convex hull. When $\epsilon>0$, $\tilde T_j^k \omega'^*_{\epsilon}$ overlaps both $\tilde T_j^{k-1} \omega'^*_{\epsilon}$ and $\tilde T_j^{k+1} \omega'^*_{\epsilon}$.}\label{figure-convh}
\end{figure}
\end{proof}
Using the Caucly estimates \ref{cauchy}, we may solve (with estimates) the (resp. ``inverse") vertical cohomological operator defined as follows.
\begin{defn}\label{def-cohom-op}
We define the (resp. ``inverse") vertical cohomological operator on $\tilde{\cA}_{\epsilon,r}^d$ (resp. $\tilde{\cA'}_{\epsilon,r}^d$)(see \rd{def-domains})
\begin{equation}\label{linear-def}L_i^v(G):=G\circ\hat\tau_i-M_i G.
\end{equation}
(resp. \begin{equation}\label{linear-inv-def}L_{-i}^v(G):=G\circ\hat \tau_i^{-1}-M_i^{-1} G.
\end{equation})	
\end{defn}
The proof of the following result is similar as \cite[Proposition 4.17]{GS}.
\begin{prop}\label{cohomo}
	Assume   $N_C$
	is   
	vertically
	Diophantine.   Fix $\epsilon_0,r_0,\delta_0, \rho_0$  in $]0,1[$. Let $0<\epsilon<\epsilon_0$, $0<\rho<\rho_0$, $0<r<r_0$, $0<\delta<\delta_0$, and $\frac{\delta}{\kappa}<\epsilon$. Suppose that $F_i\in { \cA}_{\epsilon, r}$, $i=1,\ldots, n$,
	satisfy
	\begin{equation}
		\label{almost-com}   L_i^v(F_j)- L_j^v(F_i)=0\quad
	\end{equation}
	(resp. 		\begin{equation}
		\label{almost-com2}   L_{-i}^v(F_j)- L_{-j}^v(F_i)=0\quad
	\end{equation})
	defined by $(\ref{linear-def})$ (resp. (\ref{linear-inv-def})) on  $\Omega_{\epsilon,r} \cap \hat \tau_i^{-1}\Omega_{\epsilon,r} \cap \hat \tau_j^{-1} \Omega_{\epsilon,r}$ (resp. on  $\Omega_{\epsilon,r} \cap \hat \tau_i\Omega_{\epsilon,r} \cap\hat \tau_j \Omega_{\epsilon,r}$).
	There exist functions  $G \in  \tilde{ A}_{\epsilon
		{-{\delta}/{\kappa}},r e^{-\rho}}$
	(resp. $G' \in  \tilde \cA'_{\epsilon
		{-{\delta}/{\kappa}},r e^{-\rho}}$)
	such that
	\begin{align}\label{formal2q+1}
		L_i^v(G)&=F_i \ \text{on   } \Omega_{\epsilon  {-{\delta}/{\kappa}},r e^{-\rho}}.
	\end{align}
	(resp. 		 \begin{align}\label{formal2q+1}
		L_{-i}^v(G')&=F_i \ \text{on   } \Omega_{\epsilon  {-{\delta}/{\kappa}},r e^{-\rho}}.
	\end{align})
	Furthermore, $G$ satisfies
	\begin{align}\label{estim-sol}
		\|G\|_{\epsilon  {-{\delta}/{\kappa}},re^{-\rho}}&\leq  \max_i\|F_{i}\|_{\epsilon,r}(\frac{C_1}{\delta^{\tau+\nu}}+\frac{C_1}{\rho^{\tau+\nu}}),\\
		\label{estim-solcompo}
		\|G\circ \hat \tau_i\|_{\epsilon  {-{\delta}/{\kappa}},re^{-\rho}}&\leq  \max_i\|F_{i}\|_{\epsilon,r}(\frac{  C_1}{\delta^{\tau+\nu}}+\frac{C_1}{\rho^{\tau+\nu}}).
	\end{align}
	(resp.  $G'$ satisfies
	\begin{align}\label{estim-sol2}
		\|G'\|_{\epsilon  {-{\delta}/{\kappa}},re^{-\rho}}&\leq  \max_i\|F_{i}\|_{\epsilon,r}(\frac{C_1}{\delta^{\tau+\nu}}+\frac{C_1}{\rho^{\tau+\nu}}),\\
		\label{estim-solcompo2}
		\|G'\circ \hat \tau_i^{-1}\|_{\epsilon  {-{\delta}/{\kappa}},re^{-\rho}}&\leq  \max_i\|F_{i}\|_{\epsilon,r}(\frac{  C_1}{\delta^{\tau+\nu}}+\frac{C_1}{\rho^{\tau+\nu}}).
	\end{align})
	for some constant $\kappa, C_1$  that are independent of $F,\rho,\delta, r,\epsilon$ and $\nu$ that depends only on $n$ and $d$.
	
	Moreover, the solution $G$ (resp; $G'$) is unique.
\end{prop}
\begin{proof}
	We only give the proof in the vertical cohomological operator case $(\ref{linear-def})$.
	The proof of another case is similar.
	
	Since $F_i\in {\cA}_{\epsilon, r}$, we can write
	$$
	F_i(h,v)=\sum_{Q\in \mathbb{N}^d, |Q|\geq2}\sum_{P\in\Z^n}F_{i,Q,P}h^Pv^Q,
	$$
	which converges normally   for $(h,v)\in\Om_{\epsilon,r}$.
	Note that $F_{i,Q,P}$
	are vectors, and its $k$th component
	is denoted  by $F_{i,k,Q,P}$.
	For each $(Q,P)\in \mathbb{N}^d\times \Z^n$,   each $i=1,\ldots, n$, and each $j=1,\ldots ,d$, let $ i_v:=i_v(Q,P,j)$ be in $ \{1,\ldots, n\}$ such that the maximum is realized in Definition \ref{dioph}.
	Let us set
	\begin{align}
		G_j&:=\sum_{Q\in \mathbb{N}^d, 2\leq |Q|}\sum_{P\in\Z}\frac{F_{i_v,j,Q,P}}{\lambda_{i_v}^P\mu_{i_v}^Q-\mu_{i_v,j}}h^Pv^Q,\quad j=1,\ldots. d\label{sol-coh-v}.
	\end{align}
	According to (\ref{almost-com}), we have
	\begin{equation}\label{compat-coeff}
		(\lambda_{i_v}^P\mu_{i_v}^Q-\lambda_{i_v,i})F_{m,i,Q,P}=(\lambda_{m}^P\mu_{m}^Q-\lambda_{m,i}) F_{i_v,i,Q,P}.
	\end{equation}
	Therefore, using (\ref{compat-coeff}), the $i$th-component of $ L^v_m(G)$ reads
	\begin{align*}
		L^v_m(G)_i = &\sum_{Q\in \mathbb{N}^d, 2\leq|Q| }\sum_{P\in \Z^n}(\lambda_{m}^P\mu_{m}^Q-\lambda_{m,i})\frac{F_{i_v,i,Q,P}}{(\lambda_{i_v}^P\mu_{i_v}^Q-\mu_{i_v,i})}h^Pv^Q\\
		= &\sum_{Q\in \mathbb{N}^d,2\leq|Q|}\sum_{P\in \Z^n}F_{m,i,Q,P}h^Pv^Q.
	\end{align*}
	Thus we have obtained, the formal equality~:
	\begin{equation}\label{formal-lin}
		L^v_m(G) = F_m 
		,\quad m=1,\ldots, n.
	\end{equation}
	Let us estimate these solutions.
	Without loss of generality, we may assume that $\tau \geq 1$.
	According to Definition \ref{dioph} and formula (\ref{sol-coh-v}), we have
	\begin{equation}
		\max_{j}\left (|G_{j,Q,P}|\right )\leq \max_i|F_{i,Q,P}|\frac{(|P|+|Q|)^{\tau}}{D}.
	\end{equation}
	Let $(h,v)\in \Omega_{\epsilon  {-{\delta}/{\kappa}},re^{-\rho}}$. According to (\ref{dv}) and \rrem{rem-dom}, 
	we have by convexity
	\begin{align*}
		\|G_{j,Q,P}h^Pv^Q\|&\leq \max_i\|F_{i}\|_{\epsilon,r}e^{-\delta|P|- \rho |Q|}\frac{(|P|+|Q|)^{\tau}}{D}\\
		&\leq \max_i\|F_{i}\|_{\epsilon,r}e^{-\delta|P|- \rho |Q|}\frac{(|P|^\tau+|Q|^\tau)2^{\tau}}{D}\\
		&\leq \max_i\|F_{i}\|_{\epsilon,r}(e^{-\delta/2 |P|}\frac{(4\tau e)^{\tau}}{D\delta^{\tau}}+e^{- \rho/2 |Q|}\frac{(4\tau e)^{\tau}}{D\rho^{\tau}}).
	\end{align*}
	Summing over $P$ and $Q$, we obtain
	$$
	\|G\|_{\epsilon  {-{\delta}/{\kappa}},re^{-\rho}}\leq \max_i\|F_{i}\|_{\epsilon,r}(\frac{C'}{\delta^{\tau+\nu}}+\frac{C'}{\rho^{\tau+\nu}}),
	$$
	for some constants $C',\nu$ that are independent of $F,\epsilon,\delta, \rho$. Hence, $G\in
	\cA_{\epsilon  {-{\delta}/{\kappa}},re^{-\rho}}$. 
	
	Let us prove (\ref{estim-solcompo}).   Let $B:=2\max_{\ell,j}|\mu_{\ell,j}|$. Then, there is a constant $D'$ such that
	\begin{equation}\label{sd-ehanced}
		\max_{\ell\in \{1,\ldots, n\}}  \left |\lambda_{\ell}^P\mu_{\ell}^Q-\mu_{{\ell},j}\right | \geq  \frac{D'\max_k|\lambda_{k}^P\mu_{k}^Q|}{(|P|+[Q|)^{\tau}}.
	\end{equation}
	Indeed, if $\max_k|\lambda_k^P\mu_k^Q|<B$, then Definition \ref{dioph} gives (\ref{sd-ehanced}) with $D':= \frac{D}{B}$. Otherwise, if $$|\lambda_{k_0}^P\mu_{k_0}^Q|:=\max_k|\lambda_{k}^P\mu_{k}^Q|\geq B,$$
	then $|\mu_{k_0,i}|\leq \frac{B}{2}\leq \frac{|\lambda_{k_0}^P\mu_{k_0}^Q|}{2}$. Hence, we have
	$$
	\left |\lambda_{k_0}^P\mu_{k_0}^Q-\mu_{k_0,i}\right|\geq \left||\lambda_{k_0}^P\mu_{k_0}^Q|- |\mu_{k_0,i}|\right|\geq \frac{|\mu_{k_0}^P\mu_{k_0}^Q|}{2}.
	$$
	We have verified (\ref{sd-ehanced}).
	Finally, combining all cases gives us, for $m=1,\ldots, n$,
	\begin{align*}
		|[G\circ\hat \tau_m]_{QP}| &= \left| G_{Q,P}\lambda_{m}^P\mu_{m}^Q
		\right|\leq\max_{\ell}|F_{\ell,Q,P}| \frac{|\lambda_{m}^P\mu_{m}^Q|}{|\lambda_{i_v}^P\mu_{i_v}^Q-\mu_{i_v,i}|}\\
		&\leq \max_{\ell}|F_{\ell,Q,P}|\frac{|\lambda_{m}^P\mu_{m}^Q|(|P|+|Q|)^{\tau}}{D'\max_k|\lambda_{k}^P\mu_{k}^Q|}\\
		&\leq  \max_{\ell}|F_{\ell,Q,P}|\frac{(|P|+|Q|)^{\tau}}{D'}.\end{align*}
	Hence,  $\tilde G_m:=G\circ \hat \tau_m\in \cA_{\epsilon  {-{\delta}/{\kappa}},re^{-\rho}}$.   We can define  $\tilde G\in \tilde{\cA}_{\epsilon  {-{\delta}/{\kappa}},re^{-\rho}}$ such that $\tilde G=\tilde G_m \circ \hat \tau_m^{-1}$ on $\hat \tau_m(\Omega_{\epsilon,r})$.  We verify that $\tilde G$ extends to a single-valued holomorphic function of class $\tilde {\cA}_{\epsilon,r}$. Indeed,
	$\tilde G_i\circ \hat \tau_i^{-1}=\tilde G_j\circ \hat \tau_j^{-1}$ on $\hat \tau_i(\Omega_{\epsilon,r})\cap \hat \tau_j(\Omega_{\epsilon,r})$, since the latter is connected  and the two functions agree with $G$ on $\hat \tau_i(\Omega_{\epsilon,r})\cap \hat \tau_j(\Omega_{\epsilon,r})\cap \Omega_{\epsilon,r}$ that contains a neighborhood of $\tilde\Omega_\epsilon \times\{0\}$ in $\C^{n+d}$.
	
	The uniqueness follows from the uniqueness as formal solution.
\end{proof}
The next proposition seems similar to the previous one except that it aims at defining for each $i$ a solution on some domain.
\begin{prop}\label{str-cohomo}
Assume   $N_C$ is vertically strongly Diophantine. For each $1\leq i\leq n$, let 
$$
F_i(h,v)=\sum_{Q\in \mathbb{N}^d, |Q|\geq2}\sum_{P\in\Z^n}F_{i,Q,P}h^Pv^Q,
$$
be a formal power series in $v$, with Laurent series in $h$ as coefficients. We assume that they satisfy the formal relations, $1\leq i,j\leq n$, $L_i^v(F_j)=L_j^v(F_i)$, that is $\forall (k,Q,P)\in\{1,\ldots,n\}\times \N^d \times\Z^n$
$$
(\lambda_{i}^P\mu_{i}^Q-\lambda_{i,k})F_{j,k,Q,P}=(\lambda_{j}^P\mu_{j}^Q-\lambda_{j,k}) F_{i,k,Q,P}.
$$
Then, for each $1\leq i\leq n$, there exists a unique formal power series
\begin{equation}\label{formal'2q+1}
	G^{(i)}:= \left(\sum_{Q\in \mathbb{N}^d, |Q|\geq2}\sum_{P\in\Z^n}G_{k,Q,P}^{(i)}h^Pv^Q\right)_{1\leq k\leq n}\text{such that}\quad	L_i(G^{(i)})= F_i.
\end{equation}
Furthermore, if, for a given $1\leq i\leq n$, $F_i$ is holomorphic on  $ \hat \tau_i^{-2}(\Omega_{\epsilon,r})$ (resp. on  $ \hat \tau_i^{2}\Omega_{\epsilon,r} $),  then $G^{(i)}$ is holomorphic on  $\hat \tau_i^{-2}(\Omega_{\epsilon-\frac{\delta}{\kappa},re^{-\rho}})$ for any $0<\delta<\kappa\epsilon$ and $0<\rho$ and  satisfies
\begin{align}\label{estim'-sol}
	\|G^{(i)}\circ \hat \tau^{-2}_i \|_{\epsilon  {-{\delta}/{\kappa}},re^{-\rho}}&\leq \|F_{i}\|_{\hat \tau_i^{-2}(\Omega_{\epsilon,r})}(\frac{C_1}{\delta^{\tau+\nu}}+\frac{C_1}{\rho^{\tau+\nu}}),\\
	\label{estim'-solcompo}
	\|G^{(i)}\circ \hat \tau^{-1}_i\|_{\epsilon  {-{\delta}/{\kappa}},re^{-\rho}}&\leq \|F_{i}\|_{\hat \tau_i^{-2}(\Omega_{\epsilon,r})}(\frac{  C_1}{\delta^{\tau+\nu}}+\frac{C_1}{\rho^{\tau+\nu}}).
\end{align}
for some constant $\kappa, C_1$  that are independent of $\{F_i\}_i,\rho,\delta, r,\epsilon$ and $\nu$ that depends only on $n$ and $d$. Replacing $L_i^v$ by $L_{-i}^v$ yields, for each $1\leq i\leq n$ a unique formal power series $G^{(-i)}$ satisfying  $L_{-i}^v(G^{(-i)})=F_i$ with estimates 
\begin{equation}\label{estim'-sol2}
	\max_{k=1,2}\left(\|G^{(-i)}\circ \hat \tau_i^k\|_{\epsilon  {-{\delta}/{\kappa}},re^{-\rho}}\right)\leq (\frac{  C_1}{\delta^{\tau+\nu}}+\frac{C_1}{\rho^{\tau+\nu}})\|F_{i}\|_{\hat \tau_i^{2}(\Omega_{\epsilon,r})},
\end{equation}
as above if  $F_i$ is holomorphic on  $ \hat \tau_i^{2}(\Omega_{\epsilon,r})$.
\end{prop}
\begin{proof}
	Indeed, for each $1\leq i\leq n$, the formal solution $G^{(i)}$ for the $i$th vertical cohomological equation \re{formal'2q+1} is given by~:
	\begin{align}
		G_j^{(i)}&:=\sum_{Q\in \mathbb{N}^d, 2\leq |Q|}\sum_{P\in\Z}\frac{F_{i,j,Q,P}}{\lambda_{i}^P\mu_{i}^Q-\mu_{i,j}}h^Pv^Q,\quad j=1,\ldots. d\label{sol'-coh-v}.
	\end{align}
	We recall that under strong Diophantine condition none of the denominator vanishes.
	Thus the above formal solution is meaningful.
	Moreover, 
	Diophantine inequality yields
	\begin{equation}
		|G_{i,Q,P}|\leq |F_{i,Q,P}|\frac{(|P|+|Q|)^{\tau}}{D}.
	\end{equation}
	The rest of the proof is identical to that of \rp{cohomo}.
	\end{proof}
\begin{remark}\label{i-domain}
	Assume that $\tilde N_C$ is vertically strongly Diophantine.
	Let $F_i\in { \cA}_{\epsilon, r}$, $i=1,\ldots, n$. According to \rd{def-domains}, each $F_i$ is holomorphic in a neighborhood of $\Om_{\epsilon,r}$. Let us express $F_i$ as formal power series
	$$
	F_i(h,v)=\sum_{Q\in \mathbb{N}^d, |Q|\geq2}\sum_{P\in\Z^n}F_{i,Q,P}h^Pv^Q
	$$
	normally convergent on $\Omega_{ \epsilon,r}$ satisfying to $L^v_i(F_j)=L^v_j(F_i)$ for all $1\leq i,j\leq n$.
	Assume furthermore that, for each $1\leq i\leq n$, $F_i$ is also holomorphic a neighborhood of $\hat\tau_i^{-2}(\Om_{\epsilon,r})\cup \hat\tau_i^{-1}(\Om_{\epsilon,r})$.\\
	Let us set for $\ell\in\N$ and $1\leq i\leq n$~:
	 \begin{align}
	 	\tilde\Om^{(\pm \ell)}_{i,\epsilon,r}&:=\Om_{\epsilon,r}\cup\hat\tau_i^{\pm 1}(\Om_{\epsilon,r})\cup \cdots\cup \hat\tau_i^{\pm \ell}(\Om_{\epsilon,r})\nonumber\\
	 	\tilde\Om^{(\pm \ell)}_{\epsilon,r}&=\cup_{i=1}^n\tilde\Om^{(\pm \ell)}_{i,\epsilon,r}\label{domainsl}.
	 \end{align}
	Hence, for each $1\leq i\leq n$, $F_i$ is holomorphic in a neighborhood of $\tilde\Om^{(-2)}_{i,\epsilon,r}$.
	On the one hand, according to \rp{cohomo}, there exist a unique solution $G$ holomorphic on $\tilde\Om_{\epsilon-\frac{\delta}{\kappa},re^{-\rho}}^{(-1)}\cup\tilde\Om_{\epsilon-\frac{\delta}{\kappa},re^{-\rho}}^{(1)} $ satisfying to $L^v_i(G)=F_i$ for all $1\leq i\leq n$. On the other hand, according to \rp{str-cohomo}, for each $1\leq i\leq n$, there exists a unique solution $G^{(i)}$ holomorphic on neighborhood of $\hat\tau_i^{-2}(\Om_{\epsilon-\frac{\delta}{\kappa},re^{-\rho}})\cup \hat\tau_i^{-1}(\Om_{\epsilon-\frac{\delta}{\kappa},re^{-\rho}})$ to equation $L^v_i(G^{(i)})=F_i$. Since, for each $i$, $\hat\tau_i^{-2}(\Om_{\epsilon-\frac{\delta}{\kappa},re^{-\rho}})\cup \hat\tau_i^{-1}(\Om_{\epsilon-\frac{\delta}{\kappa},re^{-\rho}})$ intersects $\tilde\Om^{(-1)}_{\epsilon-\frac{\delta}{\kappa},re^{-\rho}}\cup \tilde\Om^{(1)}_{\epsilon-\frac{\delta}{\kappa},re^{-\rho}}$ along a single connected component, $G^{(i)}$ is the holomorphic extension of $G$ on $\hat\tau_i^{-2}(\Om_{\epsilon-\frac{\delta}{\kappa},re^{-\rho}})\cup \hat\tau_i^{-1}(\Om_{\epsilon-\frac{\delta}{\kappa},re^{-\rho}})$. Hence, $G$ is a holomorphic function on neighborhood of $\tilde\Om^{(-2)}_{\epsilon-\frac{\delta}{\kappa},re^{-\rho}}$ with estimates, for each $1\leq i\leq n$
	\begin{equation}\label{estim+i-20}
		\sup_{\tilde\Om^{(-1)}_{\epsilon-\frac{\delta}{\kappa},re^{-\rho}}\cup\tilde\Om^{(1)}_{\epsilon-\frac{\delta}{\kappa},re^{-\rho}}}\|G(h,v)\|\leq \sup_{\Om_{\epsilon,r}}\|F_i\|(\frac{  C_1}{\delta^{\tau+\nu}}+\frac{C_1}{\rho^{\tau+\nu}}).
	\end{equation}
	\begin{equation}\label{estim+i-2}
		\sup_{\hat\tau_i^{-2}(\Om_{\epsilon-\frac{\delta}{\kappa},re^{-\rho}})\cup \hat\tau_i^{-1}(\Om_{\epsilon-\frac{\delta}{\kappa},re^{-\rho}})}\|G\|\leq \sup_{ \hat\tau_i^{-2}(\Om_{\epsilon,r})}\|F_i\|(\frac{  C_1}{\delta^{\tau+\nu}}+\frac{C_1}{\rho^{\tau+\nu}}).
	\end{equation}
	Similarly, assume that, for each $1\leq i\leq n$, $\tilde F_i$ is holomorphic in a neighborhood of $\tilde\Om^{(2)}_{i,\epsilon,r}$ and satisfying to $L^v_{-i}(\tilde F_j)= L^v_{-j}(\tilde F_i)$ for all $i,j$. Then there exists a unique $\tilde G$ holomorphic on neighborhood of $\tilde\Om^{(2)}_{\epsilon-\frac{\delta}{\kappa},re^{-\rho}}$ satisfying to $L^v_{-i}(\tilde G)=\tilde F_i$ with estimates similar to the above ones and written compactly as~: for each $1\leq i\leq n$
	\begin{equation}\label{estim-i+2}
		\sup_{(h,v)\in\tilde\Om^{(2)}_{i,\epsilon-\frac{\delta}{\kappa},re^{-\rho}}}\|\tilde G(h,v)\|\leq \sup_{\Om_{\epsilon,r}\cup \hat\tau_i^{2}(\Om_{\epsilon,r})}\|\tilde F_i\|(\frac{  C_1}{\delta^{\tau+\nu}}+\frac{C_1}{\rho^{\tau+\nu}}).
	\end{equation}
	\end{remark}
\section{Proof of the main result}

We are interested in the existence of a non-singular
holomorphic foliation of the germ of neighborhood of $C$ in $M$ having $C$
as a compact leaf. We refer to it as a``horizontal foliation” if exists.
The above-mentioned ``horizontal foliation” will be obtained if
we can find $\Phi= \mathrm{Id}+\phi$ be a biholomorphism of $\Omega_{\epsilon_1, r_1}$ (to be chosen) such that for any $i$,
\begin{equation}
\Phi \circ \tilde \tau_i =\tau_i \circ \Phi \label{Phi}
\end{equation}
for some biholomorphism of  $\Omega_{\epsilon_1, r_1}$ (to be chosen)
$$ \tilde \tau_i(h,v)=(\tilde \tau_i^h(h,v), M_iv)$$
such that $(M,C)$ is biholomorphic to the quotient of $\Omega_{\epsilon_1,r_1}$ by  $\tilde\tau_1,\dots,\tilde \tau_n$.
Such $ \tilde \tau_i$ is called a vertical linearization of $(M,C)$.  In fact, the codimension $d$ ``horizontal foliation” can be defined $v=$constant.
Note that $C$ is the leaf defined 
by $\{v=0\}$.

\begin{defn}
	A germ of neighborhood $(M,C)$ is vertically linearized up to order $m$ if for all $1\leq i\leq n$ 
	$$
	\tau_i(h,v)=(T_ih+\tau_{i,\geq 2}^{*,h}, M_iv+\tau_{i,\geq m}^{*,v})
	$$ 
	where $\tau_{i,\geq k}^{*,\bullet}$ denotes an analytic function on a neighborhood of $\Omega_{\epsilon,r}$, for some $0<\epsilon, r$ , vanishing at $v=0$ at order $\geq k$.
\end{defn}
\begin{prop}\label{linear-all-order}
	Under the Diophantine condition \re{dv}, a germ of neighborhood $(M,C)$ is vertically linearizable up to order $m$,  for any $m\geq 2$.
\end{prop}
\begin{proof}
We argue by induction on $m$.
	Assume that a germ of neighborhood $(M,C)$ is vertically linearized up to order $m\geq 2$, that is all $1\leq i\leq n$ 
	$\tau_i(h,v)=(T_ih+\tau_{i}^{*,h}, M_iv+\tau_{i}^{*,v})$ with $\text{ord}_{v=0}\tau_{i}^{*,v}\geq m$ and $\tau_i^*=(\tau_{i}^{*,h},\tau_{i}^{*,v})\in \cA^{n+d}_{\epsilon,r}$. Let us show that 
	\begin{equation}\label{compat}
	L_i^v([\tau_{j}^{*,v}]_{m})= L_j^v([\tau_{i}^{*,v}]_{m})
	\end{equation} on $\Om_{\epsilon',r'}^{i,j}=\Om_{\epsilon',r'}\cap \tau_i^{-1}(\Om_{\epsilon',r'})\cap\tau_j^{-1}(\Om_{\epsilon',r'})$ for all $1\leq i,j\leq n$ for some $0<\epsilon'<\epsilon$, $0<r'<r$. Indeed, recalling that $\tau_i\circ\tau_j= \tau_j\circ\tau_i$ on a neighborhood  $( \Omega_{\epsilon} \cap \hat \tau_i^{-1}\Omega_{\epsilon} \cap \hat \tau_j^{-1} \Omega_{\epsilon}) \times \{0\}$, the vertical component of  which reads
	\begin{eqnarray*}
		M_i\tau_j^{*,v}-\tau_j^{*,v}(\hat \tau_i)&=& 	M_j\tau_i^{*,v}-\tau_i^{*,v}(\hat \tau_j)\\
		&&+ \left(\tau_i^{*,v}(\hat\tau_j+\tau_j^{*})- \tau_i^{*,v}(\hat\tau_j)\right)\\
			&&- \left(\tau_j^{*,v}(\hat\tau_i+\tau_j^{*})- \tau_j^{*,v}(\hat\tau_i)\right).
	\end{eqnarray*}
The Taylor expansion at $v=0$ of one of the two last line is
$$
D_h\tau_i^{*,v}(\hat\tau_j)\tau_j^{*,h}+D_v\tau_i^{*,v}(\hat\tau_j)\tau_j^{*,v}+\text{h.o.t}.
$$
The first term is order $\geq m+2$, while the second is order $\geq 2m-1$. 
Hence, the last two lines are of order $\geq m+1$ at $v=0$ so that the truncation at degree $m$ of the equality gives the result.
According to \rp{cohomo}, there exists a solution $G\in \tilde{\cA}_{\epsilon-\frac{\delta}{\kappa},re^{-\rho}}^d$ to 
$$
L_i^v(G)=-[\tau_{i}^{*,v}]_{m}, \quad i=1,\ldots, n.
$$
Furthermore, since family $\{\hat\tau_i\}_i$ is non-resonant, the solution $G$ is unique and homogeneous of degree $m$ in the vertical direction.
Let us set $\Phi(h,v):=(h,v+G(h,v))$. Then $\tilde \tau_i:=\Phi\tau_i\Phi^{-1}$ is vertically linearized up to order $m+1$ on an appropriate domain. Indeed, we have
\begin{eqnarray*}
	M_iv+M_iG(h,v)+\tilde\tau_i^{*,v}(h,v+G(h,v))&= & (M_iv+\tau_i^{*,v})+G(\hat\tau_i+\tau_i^{*})\\
M_iG(h,v)+\tilde\tau_i^{*,v}(h,v)+(\tilde\tau_i^{*,v}(h,v+G)-\tilde\tau_i^{*,v}(h,v))	&=&\tau_i^{*,v}(h,v)+ G(\hat\tau_i)\\
&&+(G(\hat\tau_i+\tau_i^{*})-G(\hat\tau_i)).
\end{eqnarray*}
Hence we have on an appropriate domain
\begin{align*}
\tilde\tau_i^{*,v}(h,v) =& L_i^v(G)+	[\tau_i^{*,v}]_m\\
&+ (G(\hat\tau_i+\tau_i^{*})-G(\hat\tau_i)) +(\tau_i^{*,v}-[\tau_i^{*,v}]_m)\\
&-(\tilde\tau_i^{*,v}(h,v+G(h,v))-\tilde\tau_i^{*,v}(h,v)).
\end{align*}
According to \cite[Lemma 4.18]{GS}, the appropriate domain on which the previous equality holds is of the form $\Omega_{\tilde \epsilon,\tilde r}$ for some $0<\tilde \epsilon<\epsilon$ and $0<\tilde r<r$ if $\tau^{*,v}$ is small enough on $\Omega_{ \epsilon,r}$. In particular, it is a product domain with a neighborhood of $0\in \C^d$.
As the first line of the right hand side is zero by construction, we check that the other two lines are of order $\geq m+1$. Thus, $\tilde\tau_i^{*,v}$ is of order $\geq m+1$ at $v=0$.
\end{proof}
By the identity theorem, the equation $(\ref{Phi})$ implies that
\begin{equation}
\Phi \circ \tilde \tau_i^{-1} =\tau_i^{-1} \circ \Phi \label{Phi_inv}
\end{equation}
whenever both sides are well defined.



We want to find $\Phi(h,v)=(h,v)+\phi(h,v)$ with $\phi$ of the form $\phi(h,v)=(0, \phi^v(h,v))$.
Denote
$$\tau_i^{\pm 1}=(\tau^h_{i,\pm}, \tau^v_{i,\pm}).$$
Assume that they are defined on $\Omega_{\epsilon_0, r_0}$
for suitable choice (sufficiently small) $\epsilon_0, r_0$ such that $(C,M)$ is biholomorphic to the quotient of $\Omega_{\epsilon_0, r_0}$ by $\tau_j(1 \leq j \leq n)$.
Assume also that they are all defined on (see \re{domainsl})
\begin{equation}\label{domains-def-init}
 \tilde\Omega^{(-2)}_{\epsilon_0, r_0}\cup\tilde\Omega^{(2)}_{\epsilon_0, r_0}
\end{equation}
for the same choice $\epsilon_0, r_0$.

Using the condition that $N_C$ is unitary,
we have that
$\hat\tau_i(\Omega_{\epsilon_0, r_0})=T_i \Omega_{\epsilon_0} \times \Delta_{r_0}^d $.

Applying \rl{convh} with  $\epsilon_0$, there exists $\eta>0$ depending on $\epsilon_0$ such that $$ \cup_{i=0}^n T_i \Omega_{\epsilon_0+\eta} \cup \cup_{i=0}^n T_i^{-1} \Omega_{\epsilon_0+\eta}$$
is contained in the preimage of the convex hull $\mathrm{Conv}(\omega'^*_{\epsilon_0})$ of $\omega'^*_{\epsilon_0}$ under $\varphi$ with the same notations of \rl{convh}.

Define the higher order perturbations
$$\tilde \tau_{i, \pm}^{*,h}:=\tilde \tau_{i,\pm}^{h}-T_i^{\pm 1}, \tau_{i,\pm}^{*,h}:= \tau_{i, \pm}^{h}-T_i^{\pm 1},$$
$$\tau_{i, \pm}^{*,v}:= \tau_{i, \pm}^{v}-M_i^{\pm 1}.$$
The horizontal part of equation $(\ref{Phi})$ is given by
\begin{equation}
T_i^{\pm 1} h+\tilde \tau_{i, \pm}^{*,h}(h,v)=T_i^{\pm 1} h+\tau_{i, \pm}^{*,h}(h,v+\phi^{v}(h,v)),
\end{equation}
that is
\begin{equation}\label{horizontal}
\tilde \tau_{i, \pm}^{*,h}(h,v)=\tau_{i, \pm}^{*,h}(h,v+\phi^{v}(h,v)),
\end{equation}
The vertical part of equation $(\ref{Phi})$ is given by
\begin{equation}
M_i^{\pm 1} v+\phi^v(T_i^{\pm 1} h+\tilde \tau_{i, \pm}^{*,h}(h,v),M_i^{\pm 1} v) =M_i^{\pm 1} (v+\phi^v(h,v))+\tau_{i, \pm}^{*,v}(h,v+\phi^{v}(h,v)),
\end{equation}
that is
\begin{equation}\label{vertical}
\phi^v(T_i^{\pm 1} h+\tilde \tau_{i, \pm}^{*,h}(h,v),M_i^{\pm 1} v) =M_i^{\pm 1} \phi^v(h,v)+\tau_{i, \pm}^{*,v}(h,v+\phi^{v}(h,v)).
\end{equation}
We recall from \rd{def-cohom-op}, the (resp. ``inverse") vertical cohomological operator~: $L_i^v(G):=G(T_i h, M_i v)-M_i G(h,v)$ (resp. $L_{-i}^v(G):=G(T_i^{-1} h, M_i^{-1} v)-M_i^{-1} G(h,v)$).

Using equations $(\ref{horizontal})$ and $(\ref{vertical})$, we have
\begin{align}\label{linear+}
L_i^v(\phi^v)(h,v)&=\tau_{i, +}^{*,v}(h,v+\phi^v(h,v))\nonumber\\ &-\left(\phi^v(T_ih+\tau_{i,+}^{*,h}(h, v+\phi^v(h,v)), M_i v)-\phi^v(T_ih, M_i v)\right).
\end{align}
This is defined by developing the horizontal and vertical parts of equation (\ref{Phi_inv}).
\begin{remark}
Under the Diophantine assumption, the formal linearization $\Phi$ is unique. Indeed, if $\Phi$ is a equivalence between two vertically linear neighborhood, then one has $\tau_{i,+}^{*,v}=0$ in the previous equality (\ref{linear+}). If $\phi^v$ is of order $m\geq 2$ at $v=0$ then $\phi^v(T_ih+\tau_{i,+}^{*,h}(h, v+\phi^v(h,v)), M_i v)-\phi^v(T_ih, M_i v)$ is of order $\geq m+2$. Then, the previous equation (\ref{linear+}) give $L_i^v([\phi^v]_m+[\phi^v]_{m+1})=0$ for all $1\leq i \leq n$. Hence $[\phi^v]_m=[\phi^v]_{m+1}=0$ and $\phi$ is of order $\geq m+2$.
An induction on $m$ demonstrates the uniqueness of the formal linearization.
\end{remark}


Using equations $(\ref{horizontal})$ and $(\ref{vertical})$, we have
\begin{align}\label{linear-}
L_{-i}^v(\phi^v)(h,v)&=\tau_{i, -}^{*,v}(h,v+\phi^v(h,v))\\ &-(\phi^v(T_i^{-1}h+\tau_{i,-}^{*,h}(h, v+\phi^v(h,v)), M_i^{-1} v)-\phi^v(T_i^{-1}h, M_i^{-1} v)).\nonumber
\end{align}

 According to \rp{linear-all-order}, there exists a formal (power series in $v$, with holomorphic coefficients in some $\Omega_{\epsilon'}$) solution to both \re{linear+} and \re{linear-}.
In the following, we will estimate the $L^\infty$ norm to show that this formal solution is in fact convergent.
To do so, we will follow the majorant method in \cite[Section 3.3]{GS21}.
Denote
\begin{equation}\label{(I)}
(I)_{i,\pm}:=\tau_{i, \pm}^{*,v}(h,v+\phi^v(h,v)) ;
\end{equation}
\begin{equation}\label{(II)}
(II)_{i,\pm}:=\phi^v(T_i^{\pm 1}h+\tau_{i, \pm}^{*,h}(h, v+\phi^v(h,v)), M_i^{\pm 1} v)-\phi^v(T_i^{\pm 1}h, M_i^{\pm 1} v).
\end{equation}
The major difference compared to \cite[Section 3.3]{GS21} will be the estimates for $(II)$.

In the following, we will estimate $[\phi^v]_k(k \geq 2)$ by induction on $k$ (which gives the estimate for $\phi^v$.)
By identity theorem, we get the same $\phi^v$ either by the vertical cohomological operator or the inverse vertical cohomological operator if the solutions are holomorphic.

Let $0<r_1$ and $0<\epsilon_1$ be positive constants to be chosen below sufficiently small. Let us define the sequences  $r_{m+1}=r_me^{-\frac{1}{2^{m}}}$ and $\epsilon_{m+1}=\epsilon_m- \epsilon_1\frac{\eta}{ 2^{m}\kappa}$ for $m>1$ and some $\eta<\frac{\kappa}{2}$ sufficiently small. We have  $r_{m+1}:=r_1e^{-\sum_{k=1}^m\frac{1}{2^k}}$ and $\epsilon_{m+1}:= \epsilon_1-\epsilon_1\sum_{k=1}^m\frac{\eta}{2^m  \kappa}$ for $m\geq 1$. We have, for $m\geq 1$
	\begin{equation}\label{minorants}
		r_m>r_1e^{-1}=:r_{\infty},\quad  \epsilon_{m}>\epsilon_1(1-\frac{\eta}{\kappa})=:\epsilon_{\infty}>\frac{\epsilon_1}{2}.
	\end{equation} 
 Let us choose the value $\frac{\eta}{\kappa}<\frac{1}{2}$ to be the smallest value $\tilde \eta$ from \rl{cauchy-transl} corresponding to $\frac{\epsilon_1}{2}<\tilde\epsilon<\epsilon_1$.

Our goal to find germs of holomorphic function at 0
$$A(t)=\sum_{k \geq 2} A_k t^k,$$
and for $1 \leq i \leq n$,
$$B^{\pm e_i}(t)=\sum_{k \geq 2} B^{\pm e_i}_k t^k, 1 \leq i \leq n,$$
such that
\begin{equation}\label{goal}
\sup_{(h,v) \in\cup_{i=0}^n \hat\tau_i^{\pm 1}(\Omega_{\epsilon_k,r_k}) }| [\phi^v]_k(h,v)|  \leq A_k \eta_k,
\end{equation}
\begin{equation}\label{goal'}
\sup_{(h,v) \in\hat\tau_i^{\pm 1}(\Omega_{\epsilon_k,r_k}) \cup \hat\tau_i^{\pm 2}( \Omega_{\epsilon_k,r_k})} | [\phi^v]_k(h,v)| \leq B^{\pm e_i}_k \eta_k,
\end{equation}
for suitable chosen $ r_1$ (sufficiently small).
Here, the sequence $\{\eta_m\}_{m\geq 1}$ is defined by $\eta_1=1$ and for $m\geq 2$
\begin{equation}
	\eta_m
	:=\frac{C_1}{\eta^{\tau+\nu}} 2^{m(\tau+\nu)}
	\max_{m_1+\cdots +m_p+s=m} \eta_{m_1}\cdots \eta_{m_p}.\label{def-eta}
	\end{equation}
 where the constant $C_1$ is defined in Proposition \ref{cohomo}. Here we have $1\leq m_i\leq m-1$ and $s\in \N$.
We have 
\begin{equation}
\eta_m\leq\max_{1\leq s\leq m-1}\left(\frac{C_1}{\eta^{\tau+\nu}}\right)^{m-s+1}  2^{(\tau+\mu)(2m-s)}
\leq D^m,
\end{equation}
for some positive constant $D$.

Define as formal series
$$J^{m-1}A(t):=A_2t^2+\cdots+ A_{m-1}t^{m-1},$$
$$A(t)=\sum_{m \geq 2} A_m t^m,\quad B^{\pm e_i}(t)=\sum_{m \geq 2} B^{\pm e_i}_m t^m.
$$
$$$$
The idea is the following. Consider the Taylor development
$$\phi^v(h,v)=\sum_{Q \in \N^d, |Q| \geq 2} \phi_Q(h)v^Q.$$
Let $[\phi^v]_k$ be the homogenous degree $k$ part of $\phi^v$
$$[\phi^v]_k(h,v)=\sum_{Q \in \N^d, |Q| =k} \phi_Q(h)v^Q.$$
In  the following, we will alway denote $[\bullet]_k$ to indicate the homogenous degree $k$ part of some serie in $v$.
Notice that
$$[L_{\pm i}^v(\phi^v)]_k=L_{\pm i}^v [\phi^v]_k.$$
The degree 2 part  is
\begin{equation}\label{equ2+}
L_i^v([\phi^v]_2)= -[\tau_{i,+}^{*,v}(h,v)]_2,\end{equation}
whose right-handed-side term is independent of $\phi^v$.
Similarly, the degree 2 part for the inverse vertical cohomological operator is
\begin{equation}\label{equ2-}
	L_{-i}^v([\phi^v]_2)= -[\tau_{i,-}^{*,v}(h,v)]_2=M_i^{-1}[\tau_{i,+}^{*,v}(h,v)]_2\circ\hat\tau_i^{-1},
\end{equation}
whose right-handed-side term is independent of $\phi^v$. According to \re{compat} and \rp{cohomo}, these two sets of equations on $\Omega_{ \epsilon_2,r_2}$ have the same unique solution $[\phi^v]_2$ on $\tilde\Omega^{(-1)}_{ \epsilon_2,r_2}\cup \tilde\Omega^{(1)}_{ \epsilon_2,r_2}=\cup_{i=1}^n\cup_{k=-1}^1\hat\tau_i^k(\Omega_{ \epsilon_2,r_2})$, bounded there by  $2C_1\left(\frac{2}{\eta\epsilon_1}\right)^{\tau+\nu}\max_i\|[\tau_{i,+}^{*,v}(h,v)]_2\|_{\epsilon_1,r_1}$(see e.g. \re{estim-sol}).
 
For each $1\leq i\leq n$, let us obtain a bound of $G$ on $\hat\tau_i^2(\Omega_{ \epsilon_2,r_2})$ (resp. $\hat\tau_i^{-2}(\Omega_{ \epsilon_2,r_2})$). Considering equation \re{equ2-} (resp. \re{equ2+}) on $\hat\tau_i^2(\Omega_{ \epsilon_2,r_2})$ (as the right hand side is well defined according to \re{domains-def-init}) (resp. $\hat\tau_i^{-2}(\Omega_{ \epsilon_2,r_2})$), \rp{str-cohomo} and \rrem{i-domain} provide a solution which analytically continued $G$ on $\hat\tau_i^2(\Omega_{ \epsilon_2,r_2})\cup\hat\tau_i(\Omega_{ \epsilon_2,r_2})$ (resp.  $\hat\tau_i^{-2}(\Omega_{ \epsilon_2,r_2})\cup\hat\tau_i^{-1}(\Omega_{ \epsilon_2,r_2})$) bounded there by $2C_1\left(\frac{2}{\eta\epsilon_1}\right)^{\tau+\nu}\|[\tau_{i,-}^{*,v}]_2\|_{\hat{\tau}^2_i(\Omega_{\epsilon_1,r_1})}$ (resp. $2C_1\left(\frac{2}{\eta\epsilon_1}\right)^{\tau+\nu}\|[\tau_{i,+}^{*,v}]_2\|_{\hat{\tau}^{-2}_i(\Omega_{\epsilon_1,r_1})}$). We set
\begin{align*}
	A_2 &:= \max_i\|[\tau_{i,+}^{*,v}(h,v)]_2\|_{\epsilon_1,r_1}\\
	 B^{\pm e_i}_2 &:= \|[\tau_{i,\mp}^{*,v}(h,v)]_2\|_{\hat{\tau}_i^{\pm 2}(\Omega_{\epsilon_1,r_1})}.
	\end{align*}
We proceed by induction on $m\geq 2$ as 
	Taylor expansion at $v=0$ of \re{linear+} shows that for any $m \geq 2$,
	\begin{equation}\label{degm}L_i^v [\phi^v]_m=P_i(h; v, [\phi^v]_2, \cdots, [\phi^v]_{m-1})\end{equation}
	where $P_i(h; v, [\phi^v]_2, \cdots, [\phi^v]_{m-1})$ is analytic in $h\in \Omega_{ \epsilon_{m-1},r_{m-1}}$ and polynomial in $v$, $[\phi^v]_2$, ..., $[\phi^v]_{m-1}$.
To obtain the estimate (\ref{goal}) of homogeneous part of degree $m$, $[\phi^v]_m$, on $\tilde\Omega^{(-1)}_{ \epsilon_{m},r_{m}}\cup \tilde\Omega^{(1)}_{ \epsilon_{m},r_{m}}$, we invert and estimate the common solution $[\phi^v]_m$ to all vertical cohomological operators $L_i^v$ (resp. inverse vertical cohomological operators $L_{-i}^v$) from equation \re{degm}, $i=1,\ldots,n$. To do so, it is sufficient by Proposition \ref{cohomo},
to estimate the norm  of the homogeneous part of degree $m$ of its right hand side $(I)_{i,\pm}+(II)_{i,\pm}$, on $\Omega_{\epsilon_{m-1},r_{m-1}}$, for each $1\leq i\leq n$.
We notice that, according to Cauchy estimates, the $L^{\infty}$-estimate of  term $(II)_{i,\pm}$ needs estimate of terms of degree $m_1 \leq m-1$ on
$\hat\tau_i(\Omega_{\epsilon_{m_1}+ \eta,r_{m_1}})$
(resp. $\hat\tau_i^{-1}(\Omega_{\epsilon_{m_1}+ \eta,r_{m_1}})$). The latter is obtained by induction. Indeed by \rl{convh} and \rl{cauchy-transl} together with the choice of $\eta$, this domain is contained in the convex hull of the union, over $j$, of the union of 
$\hat\tau_j(\Omega_{\epsilon_{m_1},r_{m_1}})\cup \hat\tau_j^2(\Omega_{\epsilon_{m_1},r_{m_1}})$
(related to the coefficient of $B^{e_j}(t)$),
$\hat\tau_j^{-1}(\Omega_{\epsilon_{m_1},r_{m_1}}) \cup \hat\tau_j^{-2}(\Omega_{\epsilon_{m_1},r_{m_1}})$
(related to the coefficient of $B^{-e_j}(t)$)
and $\Omega_{\epsilon_{m_1},r_{m_1}} \cup \left(\cup_{k=1}^n\hat\tau_k(\Omega_{\epsilon_{m_1},r_{m_1}})\right)\cup \left(\cup_{k=1}^n\hat\tau_k^{-1}(\Omega_{\epsilon_{m_1},r_{m_1}})\right)$
(related to the coefficient of $A(t)$). The estimate of the former is thus obtained from the estimates on the latter as we remark, according to \rl{convh} that the $L^\infty-$norm on the convex hull of the union is equal to the $L^\infty-$norm on the union. The distance from $\hat\tau_i(\Omega_{\epsilon_{m_1}+ \eta,r_{m_1}})$ to the boundary of the convex hull is bounded away from $0$ independently of $m_1$ if $\eta$ is small enough.

We emphasize that the unitary flatness of the normal bundle assumption allows not to change the radius $r_{m_1}$ in the ``vertical direction" in this argument. We remark also that the usage of $B^{\pm e_i}(t)$ is necessary since the domain $\hat\tau_i(\Omega_{\epsilon_m+ \eta',r_{m}} ) $ is not contained in the convex hull of the union over $i$ of the domains $\Omega_{\epsilon_{m_1},r_{m_1}}\cup \hat\tau_i( \Omega_{\epsilon_{m_1},r_{m_1}})$ for any $m \geq m_1>0$ with some fixed $\eta'>0$ independent of $m_1,m$, as shown in Figure \ref{figure-convh}.
In particular, we have no $L^\infty-$estimate on this larger domain (i.e. $\hat\tau_i(\Omega_{\epsilon_m+ \eta',r_{m} })$ with some $\eta'>0$ independent of $m_1, m$) which is needed to apply Cauchy's estimate.

For each $1\leq i\leq n$, in order to obtain (\ref{goal'}) for a suitable $B_m^{-e_i}$ (resp. $B_m^{e_i}$), we invert and estimate the solution of the vertical cohomological operator $L_i^v$ (resp. $L_{-i}^v$) from equation \re{degm}. To do so, it is sufficient by Proposition \ref{str-cohomo},
to estimate the norm of the homogeneous part of degree $m$ of its right hand side, $(I)_{i,+}+(II)_{i,+}$ (resp. $(I)_{i,-}+(II)_{i,-}$), on $\hat\tau_i^{-2}(\Omega_{\epsilon_{m-1},r_{m-1}})$ (resp. $\hat\tau_i^{2}(\Omega_{\epsilon_{m-1},r_{m-1}})$).
We notice that, according to Cauchy estimates, the $L^{\infty}$-estimate of  term $(II)_{i,+}$ (resp. $(II)_{i,-}$) needs estimate of terms of degree $m_1 \leq m-1$ on
$\hat\tau_i^{-1}(\Omega_{\epsilon_{m_1}+ \eta,r_{m_1}})$
(resp. $\hat\tau_i(\Omega_{\epsilon_{m_1}+ \eta,r_{m_1}})$) 
which can be obtained by induction since by \rl{convh} and \rl{cauchy-transl}, this domain is contained in convex hull of the union over $j$ of the union of $\hat\tau_j(\Omega_{\epsilon_{m_1},r_{m_1}}) \cup \hat\tau_j^2(\Omega_{\epsilon_{m_1},r_{m_1}})$
(related to the coefficient of $B^{e_j}(t)$),
$\hat\tau_j^{-1}(\Omega_{\epsilon_{m_1},r_{m_1}}) \cup \hat\tau_j^{-2}(\Omega_{\epsilon_{m_1},r_{m_1}})$
(related to the coefficient of $B^{-e_j}(t)$)
and $\Omega_{\epsilon_{m_1},r_{m_1}} \cup \left(\cup_{k=1}^n\hat\tau_k(\Omega_{\epsilon_{m_1},r_{m_1}})\right)\cup \left(\cup_{k=1}^n\hat\tau_k^{-1}(\Omega_{\epsilon_{m_1},r_{m_1}})\right)$ (related to the coefficient of $A(t)$)
with estimates of norms.

 Finally, we shall show that the coefficients of $B^{\pm e_i}(t)$ and $A(t)$ of degree $m$ are bounded from above by the (non-negative) coefficient of degree $m$ of a holomorphic function of $t$, $J^{m-1}A(t)$ and $J^{m-1}B^{\pm e_i}(t)$, that is their Taylor polynomials of degree $m-1$.
This is used to proceed through a majorant method. Using the implicit function theorem for a holomorphic functional equation system of $a(t)$ and $b^{\pm e_i}(t)$, the latter being power series dominating  $A(t)$ and $B^{\pm e_i}(t)$ respectively, we can conclude that they are both holomorphic at $0$.


We will focus on the  vertical cohomological equation \re{linear+}.
The case of the inverse  vertical cohomological equation \re{linear-} is obtained similarly. We omit the "+" index in the following.
%

Let us estimate the norms of (I) and (II).

Denote $\N^d_k:=\{Q\in\N^d\colon|Q|\geq k\}$. Let $m\geq 2$, for $Q\in\mathbb N^d_2$, $|Q|\leq m$, let us set
$$
E_{Q,m}=\left\{(m_{1,1},\ldots,m_{1,q_1},\ldots,m_{d,1},\ldots,m_{d,q_d})\in\mathbb N^{|Q|}_1\colon \sum_{i=1}^d{m_{i,1}+\cdots+m_{i,q_i}=m}\right\}.
$$
For $Q\in \N^d_2$, we have
$$
\left[(v+\phi^v(h,v))^Q\right]_m= \sum_{
M\in E_{Q,m}}\prod_{j=1}^d[v_j+\phi^v_{j}]_{m_{j,1}}\cdots [v_j+\phi^v_{j}]_{m_{j,q_j}}
$$
where we have set $\phi^v=(\phi^v_{1},\ldots, \phi^v_{d})$ and $M=(m_{1,1},\ldots,m_{1,q_1},\ldots,m_{d,1},\ldots,m_{d,q_d})$.
Thus we have for any $\epsilon>0$ and $r>0$ for which $[\phi^v]_l$ is well defined $l< m$,
\begin{equation}
	\left\|\left[(v+\phi^v(h,v))^Q\right]_m\right\|_{\epsilon, r}\leq \sum_{
	M\in E_{Q,m}}\prod_{j=1}^d \left\|[v_j+\phi^v_{j}]_{m_{j,1}}\right\|_{\epsilon,r}\cdots \left\|[v_j+\phi^v_{j}]_{m_{j,q_j}}\right\|_{\epsilon,r}.
\end{equation}

Let $M'_i=(m_{1,1}^{(i)},\ldots,m_{1,q_1^{(i)}}^{(i)},\ldots,m_{d,1}^{(i)},\ldots,m_{d,q_d^{(i)}}^{(i)})\in\mathbb N^{|Q^{(i)}|}_1$
with $|Q^{(i)}|\leq m_i$
and $
m_i=\sum_{j=1}^d m_{j,1}^{(i)}+\cdots+m_{j,q_j^{(i)}}^{(i)}$, $i=1,2$.  Define the concatenation
$M'_1\sqcup M'_2$ 
to be  $(M'_1,M'_2)$.
Hence, we emphasize that the concatenation
\begin{equation}\label{concatenation}
	\left(\bigcup_{2\leq |Q_1|\leq m_1}E_{Q_1,m_1}\right)\sqcup \left(\bigcup_{2\leq |Q_2|\leq m_2}E_{Q_2,m_2}\right)\subset \bigcup_{2\leq |Q|\leq m_1+m_2}E_{Q,m_1+m_2}.
\end{equation}
By Cauchy estimate $(\ref{Reps})$ applying to $\tau_{i,j}(1 \leq j \leq n+d)$ implies that, if $\epsilon_1$ small enough, there exists $R>0$ 
such that
$$
||\tau_{i,Q,j}||_{\epsilon_1} \leq R^{|Q|}
$$
with
$$\tau_{i,j}=\sum_{Q \in \N^d} \tau_{i,Q,j} (h)v^Q$$
and $\tau_i=(\tau_{i,1}, \cdots, \tau_{i,d})$.
Without loss of generality, we may assume that same estimate holds on
$$
\tilde \Omega_{\epsilon_1}^{(-2)}\cup\tilde \Omega_{\epsilon_1}^{(2)}:=\cup_{j} T_j^{-2}\Omega_{\epsilon_1} \cup T_j^{-1} \Omega_{\epsilon_1}\cup \Omega_{\epsilon}\cup T_j^{1} \Omega_{\epsilon_1} \cup T_j^{2} \Omega_{\epsilon_1}.
$$

 We recall that $\Om_{\epsilon_j,r_j}\Subset  \Om_{\epsilon_l,r_l}$ if $l<j$. Assuming by induction that \re{goal} holds for all $m'<m$, we have
\begin{align}
	 \left\|[(I)]_m\right\|_{\epsilon_m,r_m}&\leq \sum_{|Q|=2}^{m}R^{|Q|}\sum_{M'\in E_{Q,m}}\prod_{j=1}^d
	\left\|[v_j+\phi^v_{j}]_{m_{j,1}}\right\|_{\epsilon_m,r_m}\cdots \left\|[v_j+\phi^v_{j}]_{m_{j,q_j}}\right\|_{\epsilon_m,r_m}
	\nonumber\\
	&\leq \sum_{|Q|=2}^{m}R^{|Q|}\sum_{M'\in E_{Q,m}}\prod_{j=1}^d\eta_{m_{j,1}}A_{m_{j,1}}\cdots \eta_{m_{j,q_j}} A_{m_{j,q_j}}\nonumber\\
	&\leq \left[\sum_{|Q|=2}^{m}\eta_{Q,m}R^{|Q|}(t+J^{m-1}(A(t))^{|Q|}\right]_m\nonumber\\
	&\leq  E_{m}[g_m(t)]_m,\label{estim-h'}
\end{align}
where we have set
	\begin{align*}
		\label{g(t)}
		&\eta_{Q,m}:=
		\max
		_{M'\in E_{Q,m}}\left(\prod_{i=1}^d\eta_{m_{i,1}}\cdots \eta_{m_{i,q_i}}\right),\quad E_m:=\max_{\dindice{Q\in \Bbb N^d}{2\leq |Q|\leq m}} \eta_{Q,m},\\
		& g_m(t):= \sum_{|Q|=2}^{m}R^{|Q|}(t+J^{m-1}(A(t))^{|Q|},\quad   g(t):= \sum_{|Q|\geq 2}R^{|Q|}(t+A(t))^{|Q|}.
			\end{align*}
			Here $[g(t)]_m$ denotes the coefficient of $t^m$ in the power series $g(t)$.
			We also define
		\begin{align*}
				&  G(t,U):= \sum_{|Q|\geq 2}R^{|Q|}(t+U)^{|Q|},\\
			& g^{\pm e_i}(t):= \sum_{|Q|\geq 2}R^{|Q|}(t+B^{\pm e_i}(t))^{|Q|}=G(t,B^{\pm e_i}(t)),\\
				& g_m^{\pm e_i}(t):= \sum_{|Q|=2}^{m}R^{|Q|}(t+J^{m-1}(B^{\pm e_i}(t))^{|Q|}.
			\end{align*}
We have
\begin{eqnarray*}
	[(II)]_m &=&
	\sum_{\dindice{P\in \Bbb N_1^n}{m_1+m_2=m}}\frac{1}{P!}\left[ \partial^P_h \phi^v(T_ih,M_iv)\right]_{m_1}\left[\left(\tau_{i}^{*,h}(h,v+\phi^v(h,v))\right)^P\right]_{m_2}\\
	&=&
\sum_{\dindice{P\in \Bbb N_1^n}{m_1+m_2=m}}\frac{1}{P!} (\partial^P_h \left[\phi^v\right]_{m_1})(T_ih,M_iv)\left[\left(I\right)^P\right]_{m_2}
\end{eqnarray*}
Here, both indices $m_1$ and $m_2$ are $\geq 2$  so that both $m_1$ and $m_2$ are less or equal than $m-2$. Assuming by induction that \re{goal} holds for all $m'<m$, we have
\begin{equation}\label{h'p}\nonumber
	\left\|\left[\left(I\right)^P\right]_{m_2}\right\|_{\epsilon_{m-1},r_{m-1}}\leq E_{m_2}\left[\left(\sum_{|Q|=2}^{m_2}R^{|Q|}(t+J^{m-1}(A(t))^{|Q|}\right)^{|P|}\right]_{m_2}.
\end{equation}
Indeed,
\begin{eqnarray*}
	\left[\left(I\right)^P\right]_{m_2}&=&\left[\prod_{i=1}^n((I)_{i})^{p_i}\right]_{m_2}\\
	&=&  \sum_{\sum_i(m_{i,1}+\cdots +m_{i,p_i})=m_2}\prod_{i=1}^n[(I)_{i}]_{m_{i,1}}\cdots [(I)_{i}]_{m_{i,p_i}}.
\end{eqnarray*}
Here,$(I)_{i}$ means the $i-$th component of term $(I)$.
According to $(\ref{concatenation})$ and by $(\ref{estim-h'})$, we have
\begin{eqnarray*}
 \left\|\prod_{i=1}^n[(I)_{i}]_{m_{i,1}}\cdots [(I)_{i}]_{m_{i,p_i}}\right\|_{\epsilon_{m-1},r_{m-1}}&\leq &\prod_{i=1}^nE_{m_{i,1}}\left[g_{m_{i,1}}(t)\right]_{m_{i,1}}\cdots E_{m_{i,p_i}}\left[g_{m_{i,p_i}}(t)\right]_{m_{i,p_i}}\\
	&\leq & \max_{2\leq |Q|\leq m_{2}}\eta_{Q, m_{2}}\prod_{i=1}^n\left[g_{m_{i,1}}(t)\right]_{m_{i,1}}\cdots \left[g_{m_{i,p_i}}(t)\right]_{m_{i,p_i}}.
\end{eqnarray*}
Hence, we have
$$
\sum_{\sum_i(m_{i,1}+\cdots +m_{i,p_i})=m_2} \left\|\prod_{i=1}^n[(I)_{i}]_{m_{i,1}}\cdots [(I)_{i}]_{m_{i,p_i}}\right\|_{\epsilon_{m-1},r_{m-1}}\leq  E_{m_{2}}[g(t)^{|P|}]_{m_2}.
$$

Now we estimate $\left[ \partial^P_h \phi^v(T_ih,M_iv)\right]_{m_1}$ where attention is put on the choice of domain.
By Cauchy estimate $(\ref{Reps3})$, we have  by induction on $m\geq 2$, for any $m_1 < m$, 
$$
||\left[ \partial^P_h \phi^v(T_ih,M_iv)\right]_{m_1}||_{\epsilon_{m_1}, r_{m_1}}
\leq \frac{C (C')^{|P|}(\sum_{j,\pm} B^{\pm e_j}_{m_1}\eta_{m_1}+A_{m_1}\eta_{m_1})}{C''^{\nu+|P|}}
$$
since
 $\hat\tau_i(\Omega_{\epsilon_{m_1}+\eta,r_{m_1}})$
is contained in the convex hull of the union, over $j$, of the union of 
$\hat\tau_j(\Omega_{\epsilon_{m_1},r_{m_1}})\cup \hat\tau_j^2(\Omega_{\epsilon_{m_1},r_{m_1}})$
(related to the coefficient of $B^{e_j}(t)$),
$\hat\tau_j^{-1}(\Omega_{\epsilon_{m_1},r_{m_1}}) \cup \hat\tau_j^{-2}(\Omega_{\epsilon_{m_1},r_{m_1}})$
(related to the coefficient of $B^{-e_j}(t)$)
and $\Omega_{\epsilon_{m_1},r_{m_1}} \cup\cup_{k=1}^n \hat\tau_k(\Omega_{\epsilon_{m_1},r_{m_1}})\cup\cup_{k=1}^n \hat\tau_k^{-1}(\Omega_{\epsilon_{m_1},r_{m_1}})$
(related to the coefficient of $A(t)$).
Here $C''$ is a constant independent of $m$. As a consequence, we have
\begin{align}
	||\left[\left(II\right)\right]_{m}||_{\epsilon_{m-1}, r_{m-1}}\leq &\sum_{m_1+m_2=m}\sum_{\dindice{P\in \Bbb N^n}{|P|\geq 1}}
	\frac{C (C')^{|P|}(A_{m_1}+\sum_{j,\pm} B^{\pm e_j}_{m_1})\eta_{m_1}}{C''^{\nu+|P|}}E_{m_{2}}[g(t)^{|P|}]_{m_2}\nonumber\\
	\leq & \sum_{m_1+m_2=m}\frac{C (A_{m_1}+\sum_{j,\pm} B^{\pm e_j}_{m_1})\eta_{m_1}}{C''^\nu} \left[E_{m_{2}}\sum_{\dindice{P\in \Bbb N^n}{|P|\geq 1}}
	\left(C' {g(t)}/C''\right)^{|P|}\right]_{m_2}\nonumber\\
	\leq & \frac{C }{C''^\nu} \left(\max_{m_1+m_2=m}\eta_{m_1}E_{m_2}\right)\times\nonumber\\
	&\times\left[( A(t)+\sum_{j,\pm} B^{\pm e_j}(t)) \left(\left(\frac{1}{1-C' g(t)/C''}\right)^n-1\right)\right]_m.\label{estim-h''}
\end{align}
Collecting estimates $(\ref{estim-h'})$ and $(\ref{estim-h''})$, we obtain
\begin{align*}\nonumber
		\left\|L^v_i[\phi^v]_m 	\right\|_{\epsilon_{m-1}, r_{m-1}}& \leq\left[ E_m g(t)+ \frac{C }{C''^\nu} \left(\max_{m_1+m_2=m}\eta_{m_1}E_{m_2}\right)\times\right.\\ 
	&\left.\times(A(t)+\sum_{j,\pm} B^{\pm e_j}(t)) \left(\left(\frac{1}{1-C' g(t)/C''}\right)^n-1\right)\right]_m.
\end{align*}
We solve the vertical cohomological operator for $(\ref{linear-def})$ and obtain by Proposition \ref{cohomo} the following estimate~:
\begin{align}
			\left\|[\phi^v]_m	\right\|_{\epsilon_{m}, r_{m}}&\leq \frac{ C_1}{\eta^{\tau+\nu}} 2^{m(\tau+\nu)} \left[ E_m g(t)+ \frac{C }{C''^\nu} \left(\max_{m_1+m_2=m}\eta_{m_1}E_{m_2}\right)\right.\nonumber\\
	&\left.\times (A(t)+\sum_{j,\pm} B^{\pm e_j}(t)) \left(\left(\frac{1}{1-C' g(t)/C''}\right)^n-1\right)\right]_m.\label{induction}
\end{align}
Note that we use that $\hat\tau_i\hat\tau_j=\hat\tau_j\hat\tau_i$ to apply Proposition \ref{cohomo} to $L_i^v(\phi^v)$.
Using definition $(\ref{def-eta})$, we obtain
\begin{equation}\nonumber
		\left\|[\phi^v]_m 	\right\|_{\epsilon_{m}, r_{m}}\leq \eta_m \left[  g(t)+ \frac{C }{C''^\nu}  (A(t)+\sum_{j,\pm} B^{\pm e_j}(t)) \left(\left(\frac{
		1}{1-C' g(t)/C''}\right)^n-1\right)\right]_m.
\end{equation}

For each $1\leq i\leq n$, we consider the single inverse vertical cohomological equation \re{linear-}. By \rp{str-cohomo} and \rrem{i-domain}, we obtain ,
\begin{equation}\nonumber
	\left\|[\phi^v]_m 	\right\|_{\hat\tau_i^2(\Omega_{\epsilon_{m}, r_{m}})}\leq \eta_m \left[  g^{+e_i}(t)+ \frac{C }{C''^\nu}  (A(t)+\sum_{j,\pm} B^{\pm e_j}(t)) \left(\left(\frac{
		1}{1-C' g^{+e_i}(t)/C''}\right)^n-1\right)\right]_m.
\end{equation}
By \re{estim-i+2}, we also obtain the same estimate for $	\left\|[\phi^v]_m 	\right\|_{\hat\tau_i(\Omega_{\epsilon_{m}, r_{m}})}$.
Similarly, considering, for each $1\leq i\leq n$, the single vertical cohomological equation \re{linear+}, by  \rp{str-cohomo} and \rrem{i-domain}, we also have
\begin{equation}\nonumber
		\left\|[\phi^v]_m	\right\|_{\hat\tau_i^{-2}(\Omega_{\epsilon_{m}, r_{1,m}})}\leq \eta_m \left[  g^{-e_i}(t)+ \frac{C }{C''^\nu}  (A(t)+\sum_{j,\pm} B^{\pm e_j}(t)) \left(\left(\frac{
		1}{1-C' g^{-e_i}(t)/C''}\right)^n-1\right)\right]_m.
\end{equation}
By \re{estim+i-2}, we obtain the same estimate for $	\left\|[\phi^v]_m 	\right\|_{\hat\tau_i^{-1}(\Omega_{\epsilon_{m}, r_{m}})}$.

Let us consider the functional equation system, $i=1,\ldots, n$,
\begin{equation}\label{equ-A}\nonumber
	A(t)=   G(t,A(t))+ \frac{C }{C''^\nu}  (A(t)+\sum_{j,\pm} B^{\pm e_j}(t)) \left(\left(\frac{
		1}{1-C' G(t,A(t))/C''}\right)^n-1\right),
\end{equation}
\begin{equation}\label{equ-A'}\nonumber
	B^{\pm e_i}(t)=   G(t,B^{\pm e_i}(t))+ \frac{C }{C''^\nu}  (A(t)+\sum_{j,\pm} B^{\pm e_j}(t)) \left(\left(\frac{
		1}{1-C' G(t,B^{\pm e_i}(t))/C''}\right)^n-1\right).
\end{equation}
This equation system has a unique analytic solution vanishing at the origin at order $2$ as shown by the implicit function theorem.
Notice that the coefficients of the powers of $A(t)$,  $B^{\pm e_j}(t)$ are non-negative. As $A_2=B_2^{\pm e_i}=[G(t,0)]_2> 0$, we obtain by induction that all coefficients of degree $m\geq 2$ of  $A(t)$ and  $B^{\pm e_j}(t)$ are non-negative.



We now can prove the theorem. Indeed by assumption, there are positive constants $M'',L$ such that $\eta_m\leq M''L^m$ for all $m\geq 2$. Since $A(t)$ converges at the origin, then $A_m\leq D^m$ for some positive $D$. By the majorant construction and previous estimates, we have $\left\|[\phi^v]_m 	\right\|_{\epsilon_{m}, r_{m}}\leq A_m$ for all $m\geq 2$. Hence, according to \re{minorants}, we have $$||[\phi^v]_m||_{\frac{\epsilon_1}{2},r_1e^{-1}}\leq M''(DL)^m$$ for all $m\geq 2$. Hence, $\phi^v$ converges near the torus.
This proves the theorem.

\end{document}